\documentclass[11pt]{article}
\usepackage{graphicx} \usepackage{amsmath}
\usepackage{amsfonts} \usepackage{amssymb}

\newtheorem{theorem}{Theorem}

\newtheorem{lemma}[theorem]{Lemma}

\newtheorem{proposition}[theorem]{Proposition}

\newenvironment{proof}[1][Proof]{\textbf{#1.} }{\
\rule{0.5em}{0.5em}}

\def\fidi{\hskip5pt \vrule height4pt width4pt depth0pt \par}

\begin{document}

\author{N. Ciccoli
\thanks{supported by  PRIN \emph{Azioni di gruppi su variet\'a} and GNSAGA.}
\\ Dipartimento di Matematica e Informatica\\ Universit\'a di Perugia
\and A. J.--L. Sheu \thanks{supported by the University of Kansas General
research Fund allocation \#2301 for FY 2004.}\\ Department of Mathematics\\ University of Kansas}

\title{\bf Covariant Poisson structures\\ on complex
Grassmannians}
\maketitle
\begin{abstract}
The purpose of this paper is to study covariant Poisson structures on
$Gr_k^n\mathbb C$ obtained as quotients by coisotropic subgroups of the standard
Poisson--Lie $SU(n)$. Properties of Poisson quotients allow to describe
Poisson embeddings generalizing those obtained in \cite{Sh4}.
\end{abstract}

\section{Introduction}

In \cite{Sh4} a family of covariant Poisson structures on complex
projective spaces underlying the Dijkhuizen--Noumi quantization
(\cite{DiNo}) was studied from the point of view of coisotropic
subgroups with respect to an affine Poisson structure on $SU(n)$,
providing also a description of the associated Lagrangian
subalgebras and Poisson embeddings of standard odd Poisson spheres
in non standard Poisson projective spaces.

In this paper we plan to extend those results to complex
Grassmannians (their quantum version may be found in \cite{qGr}). The
emphasis is laid even more strongly on the role played by
subgroups which are coisotropic with respect to the standard
multiplicative Poisson structure on $SU(n)$.

One reason of interest lies in the fact that every coisotropic subgroup of a
Poisson--Lie group can be quantized in such a way as to fit in a
nice duality diagram (\cite{CiGa}). Furthermore coisotropic submanifolds have recently
raised a lot of attention in the context of deformation quantization (\cite{CF,BGHW}) and played
a role in the analysis of Poisson sigma--models over group manifolds (\cite{BZ}).

In the first section we clarify the relation between coisotropic
subgroups of a Poisson--Lie group and coisotropic subgroups of
translated affine Poisson bivectors. This relates results in \cite{Sh4}
with those in the present work, allowing a natural interpretation from the foliation
point of view.

In the second section we describe the family of covariant Poisson
structures on complex Grassmannians under consideration and show how it can be obtained as
quotients by coisotropic subgroups. Such structure was first introduced in
\cite{KRR} under different methods. A specific non standard Grassmannian was studied,
recently, with Lie group methods by Foth and Lu (see \cite{Flu}).

Finally in the last section we describe a general
procedure allowing to determine Poisson embeddings of $G$--spaces.
Applying it to projective spaces we show how it recovers the whole
symplectic foliation in the standard case and the Poisson embeddings of
\cite{Sh4} in the nonstandard one.
Moving on to Grassmannians such procedure will give
embeddings of standard Poisson--Stiefel manifolds (and of other more general manifolds) in non standard
Poisson Grassmannians. In the special case of Grassmannians
$Gr_{m}^{2m}(\mathbb C)$ this will result in a Poisson embedding of the standard Poisson--Lie group $U(m)$.
Such embeddings are relevant also from the point of view of quantum spaces, where they were first identified.
It is in the study of the groupoid $C^*$--algebra ${\cal C}(\mathbb P^n_{q,c})$ carried out in \cite{Sh:cqg,Sh:qcp},
in fact, they were used to construct composition sequences for the algebra and, eventually, to compute
its $K$--groups. We expect that quantum Stiefel manifolds studied in \cite{PV} and suitable generalizations will
appear as quotients of nonstandard complex $q$--Grassmannians and allow a similar detailed analysis.

\section{Coisotropic and affine coisotropic subgroups}

\subsection{Affine Poisson structures}

Let $G$ be a given Lie group, with Lie algebra $\mathfrak g$. In the following, we use
$R_{g}$ (resp. $L_{g}$) to denote the right (resp. left)
translation action on $G$ by $g\in G$, and also all the actions
induced by it on tensors of $G$, e.g. $R_{g}\left(  v\right)
=\left(  D\left( R_{g}\right)  \right)  _{x}\left(  v\right)  \in
T_{xg}G$ for any vector $v\in T_{x}G$, and $\left(  R_{g}X\right)
\left(  x\right)  =\left(  D\left( R_{g}\right)  \right)
_{xg^{-1}}\left(  X\left(  xg^{-1}\right)  \right)  \in T_{x}G$
for any vector field $X\in\Gamma\left(  TG\right)  $, where
$D\left( R_{g}\right)  $ is the differential (a vector bundle map
on $TG$) of the diffeomorphism $R_{g}$.
Note that the right translation $R$ is an anti-homomorphism, i.e. $R_{g}%
R_{h}=R_{hg}$. Similarly we have the left translation $L_{x}$, but
$L$ is a homomorphism, i.e. $L_{g}L_{h}=L_{gh}$. Given any $2$--tensor $\rho:G\to\wedge^2 TG$
let $\tilde\rho(g):=L_g^{-1}\rho(g)\in\wedge^2\mathfrak g$.

First we recall the following facts for an (alternating) 2-tensor field $\rho$ on a Lie
group $G$ (see \cite{Lu, We:aps}).
\begin{itemize}
\item[(1)] $\rho$ is called multiplicative if \[ \rho\left(
gh\right)  =L_{g}\left(  \rho\left(  h\right)  \right)
+R_{h}\left(  \rho\left(  g\right)  \right)  . \] (Note that
$\rho\left(  e\right)  =0$ if $\rho$ is multiplicative,  where
$e$ is the unit of $G$.)
\item[(2)] $\rho$ is called affine if \[ \rho\left(  gh\right)  =L_{g}\left(
\rho\left(  h\right)  \right) +R_{h}\left(  \rho\left(  g\right)
\right)  -L_{g}R_{h}\left(  \rho\left( e\right)  \right)\, . \]
\item[(3)] $\rho$ is affine if and only if
$\pi:=\rho-\left(  \rho\left(  e\right) \right)  ^{l}$ is
multiplicative, where $X^{l}$ denotes the left-invariant tensor
field generated by $X\in\wedge^{2}\mathfrak{g}$.
(Note that for any
2-tensor field $\pi$ with $\pi(e)=0$ and $X\in\wedge^{2}\mathfrak{g}$, if
$\rho:=\pi+X^{l}$, then $\rho\left(  e\right)  =X$ and hence
$\pi=\rho-\left(  \rho\left(  e\right) \right)  ^{l}$. So all
affine $\rho$ are of the form $\rho=\pi+X^{l}$ for some
multiplicative $\pi$ and $X\in\wedge^{2}\mathfrak{g}$.)
\item[(4)] If $\pi$
is a Poisson-Lie structure on $G$, then $\rho=\pi+X^{l}$ (with
$\rho\left(  e\right)  =X$) is affine for any
$X\in\wedge^{2}\mathfrak{g}$ (but may not be Poisson);
in this case $\rho$ is also Poisson if and only if $dX=\frac 1 2 [X,X]$ (\cite{DaSo}).
\item[(5)] If $\rho$ is affine Poisson, then $\pi:=\rho-\left(
\rho\left( e\right) \right)  ^{l}$ is multiplicative \textbf{and
Poisson}. (But the converse may not be true, cf. (4) above.)
\item[(6)] \textbf{Given }$\rho$\textbf{ Poisson}, we have that $\rho$ is
affine Poisson if and only if $\pi:=\rho-\left(  \rho\left(
e\right)  \right)  ^{l}$ is multiplicative Poisson (or
Poisson-Lie).
\end{itemize}
Now we show that if $\rho$ is affine Poisson and
$\rho\left(  \sigma\right) =0$ for some point $\sigma\in G$, then
$R_{\sigma^{-1}}\rho$ is Poisson-Lie.

\begin{lemma}. If $\rho$ is affine Poisson, then $R_{\sigma}\rho$
is also affine Poisson for any $\sigma\in G$. (We don't assume
$\rho\left(\sigma\right)  =0$ in this lemma.)%
\end{lemma}

\proof Clearly the right translation of a Poisson structure on $G$
is still a Poisson structure. So $R_{\sigma}\rho$ is Poisson. Now
by the commutativity $R_{g}L_{h}=L_{h}R_{g}$ for all $g,h\in G$,
we get \[ \left[  R_{\sigma}\rho-\left(  \left(
R_{\sigma}\rho\right)  \left( e\right)  \right)  ^{l}\right]
\left(  g\right)  =R_{\sigma}\left(
\rho\left(  g\sigma^{-1}\right)  \right)  -L_{g}\left(  \left(  R_{\sigma}%
\rho\right)  \left(  e\right)  \right)
\]%
\[ =R_{\sigma}\left[  L_{g}\left(  \rho\left(  \sigma^{-1}\right)
\right)
+R_{\sigma^{-1}}\left(  \rho\left(  g\right)  \right)  -L_{g}R_{\sigma^{-1}%
}\left(  \rho\left(  e\right)  \right)  \right]  -L_{g}\left(
R_{\sigma }\left(  \rho\left(  \sigma^{-1}\right)  \right)
\right)
\]%
\[ =R_{\sigma}\left(  L_{g}\left(  \rho\left(  \sigma^{-1}\right)
\right) \right)  +\rho\left(  g\right)  -L_{g}\left(  \rho\left(
e\right)  \right) -L_{g}\left(  R_{\sigma}\left(  \rho\left(
\sigma^{-1}\right)  \right) \right)
\]%
\[ =\rho\left(  g\right)  -L_{g}\left(  \rho\left(  e\right)
\right)  =\left[ \rho-\left(  \rho\left(  e\right)  \right)
^{l}\right]  \left(  g\right) \] which shows that \[
R_{\sigma}\rho-\left(  \left(  R_{\sigma}\rho\right)  \left(
e\right)
\right)  ^{l}=\rho-\left(  \rho\left(  e\right)  \right)  ^{l}%
\] a multiplicative Poisson structure since $\rho$ is affine
Poisson. Thus $R_{\sigma}\rho$ is affine Poisson. \endproof

\begin{proposition}. If $\rho$ is affine Poisson and $\rho\left(
\sigma\right)  =0$ for some point $\sigma\in G$, then $R_{\sigma^{-1}}\rho$ is Poisson-Lie.
\end{proposition}

\proof $R_{\sigma^{-1}}\rho$ is affine Poisson with \[ \left(
R_{\sigma^{-1}}\rho\right)  \left(  e\right)
=R_{\sigma^{-1}}\left( \rho\left(  \sigma\right)  \right)  =0 \]
and hence $R_{\sigma^{-1}}\rho$ is multiplicative Poisson.%
\endproof

\subsection{Coisotropic subgroups}

In this section we will clarify the relation between affine Poisson structures
on Lie groups and coisotropic subgroups of Poisson--Lie groups, introducing the notion
of affinely coisotropic subgroup.

Recall that for a given Poisson manifold $(M,\pi_M)$ a coisotropic
submanifold is an embedded submanifold such that its defining ideal (i.e.
the ideal of smooth functions which are zero on the manifold) is a
Poisson subalgebra. For a given Poisson--Lie group $(G,\pi)$ a coisotropic
subgroup is a Lie subgroup $H$ which is also a coisotropic
submanifold. At the infinitesimal level, if $\delta=(D\tilde\pi)_e:{\mathfrak
g}\to \wedge^2{\mathfrak g}$ represents the cobracket and
$\mathfrak h$ is a Lie subalgebra of $\mathfrak g$ then $\mathfrak h$ can be
integrated to a coisotropic subgroup if and only if
$\delta({\mathfrak h})\subseteq{\mathfrak h}\wedge{\mathfrak g}$.

Let $\rho$ be an affine Poisson structure on the Lie group $G$ and
let $H$ be a closed (connected) subgroup. It is known that the
multiplicative Poisson structure on $G$ induces (or projects to) a
well-defined Poisson structure on $G/H$ when $H$ is a coisotropic subgroup.

The concept of a coisotropic subgroup
$H$ of an affine Poisson Lie group $\left(  G,\rho\right)  $ is more delicate,
and there is a fine distinction between ``a coisotropic subgroup'' and ``a
subgroup that is a coisotropic submanifold'' as discussed below.
First we note that the following conditions are equivalent:
\begin{itemize}
\item[(2)] $\rho\left( h\right)  -R_{h}\left(  \rho\left(  e\right)  \right)
\in L_{h}\left(  \mathfrak{h}\wedge\mathfrak{g}\right)  $ for all
$h\in H$;
\item[(3)] $\rho\left(  kh\right)  -R_{h}\left(  \rho\left(
k\right) \right)  \in L_{kh}\left(
\mathfrak{h}\wedge\mathfrak{g}\right) $ for all $h,k\in H$;
\item[(4)] $\left( \left(  D\tilde\rho\right) _{e}+\left[\tilde\rho\left(  e\right)
,\cdot\right] \right)  \left( \mathfrak{h}\right)
\subset\mathfrak{h}\wedge\mathfrak{g}$;
\item[(5)] $\rho\left( gh\right)
-R_{h}\left(  \rho\left(  g\right) \right) \in L_{gh}\left(
\mathfrak{h}\wedge\mathfrak{g}\right)  $ for all $h\in H$ and
$g\in G$.
\end{itemize}
Furthermore if $ad_{\mathfrak h}(\rho(e))\subseteq{\mathfrak h}\wedge{\mathfrak g}$
such conditions are equivalent to
\begin{itemize}
\item[(1)] $H$ is a $\rho$-coisotropic submanifold of
$G$, i.e. $\rho\left(  h\right) \in L_{h}\left(
\mathfrak{h}\wedge\mathfrak{g}\right)  $ for all $h\in H$;
\end{itemize}
In fact,
\[
R_{h}\left(  \rho\left(  e\right)  \right)  =L_{h}\left(  L_{h}^{-1}%
R_{h}\left(  \rho\left(  e\right)  \right)  \right)  =L_{h}\left(
\operatorname{Ad}_{h}^{-1}\left(  \rho\left(  e\right)  \right)
\right)
\]
and so (1) $\Leftrightarrow$ (2) if
$\operatorname{ad}_{\mathfrak{h}}\left( \rho\left(  e\right)  \right)
\subset\mathfrak{h}\wedge\mathfrak{g}$. Since $\rho$ is affine, we
have \[ L_{g}\left[  \rho\left(  h\right)  -R_{h}\left(  \rho\left(
e\right) \right) \right]  =\rho\left(  gh\right)  -R_{h}\left(
\rho\left( g\right)  \right) \]
for any $h\in H$ and $g\in G$, and
hence (2), (3), and (5) are clearly equivalent. From
\[
L_{h}^{-1}\left[ \rho\left(  h\right)  -R_{h}\left(  \rho\left(
e\right)
\right)  \right]  =L_{h^{-1}}\rho\left(  h\right)  -\operatorname{Ad}_{h}%
^{-1}\left(  \rho\left(  e\right)  \right) \] and $D\left[
\operatorname{Ad}_{h}^{-1}\left(  \rho\left(  e\right)  \right)
\right]  _{h=e}=-\operatorname{ad}_{\cdot}\left(  \rho\left(
e\right) \right)  $, it is not hard to see the equivalence of (3)
and (4).

Note that even if $\rho$ is multiplicative and a subgroup
$H$ is a $\rho$--coisotropic submanifold, a coset $gH$
of $H$ in general need not be a $\rho$-coisotropic submanifold of
$G$, but is ``\textbf{affinely }(or relatively)\textbf{\
}$\rho$\textbf{-coisotropic}'' in the sense of condition (5). Note
that in general, when $\rho\left(  e\right)  \neq0$, i.e. $\rho$ is
not multiplicative, both (1)$\Rightarrow$(2) and
(2)$\Rightarrow$(1) may not hold. We \textbf{define} a closed
subgroup $H$ of an (affine) Poisson Lie group $G$ to be a
$\rho$-\textbf{coisotropic subgroup\ }if each coset $gH$ with $g\in
G$ is an \textbf{affinely }$\rho$\textbf{-coisotropic} submanifold
of $G$, i.e. $\rho\left(  gh\right)  -R_{h}\left(  \rho\left(
g\right)  \right)  \in L_{gh}\left(  \mathfrak{h}\wedge\mathfrak{g}\right)
$ for all $h\in H$. So when $\rho$ is multiplicative, a closed
subgroup $H$ of $G$ is a $\rho$-coisotropic submanifold of $G$ if
and only if $H$ is a $\rho$-coisotropic subgroup of $G$.

\begin{theorem}. Let $\pi$ be a Poisson 2-tensor on a Lie group
$G$. For a closed Lie subgroup $H$ of $G$ and $\sigma\in G$, the
conjugate $\operatorname{Ad}_{\sigma}H$ of $H$ is
$\pi$-coisotropic if and only if $H$ is
$\pi_{\sigma}$-coisotropic, where $\pi_{\sigma}\left(  g\right)
:=R_{\sigma}\left(  \pi\left(  g\sigma^{-1}\right)  \right)  $ for $g\in G$.
\end{theorem}
\proof
Let $\tilde{\pi}\left(  g\right)  :=L_{g}^{-1}\left(
\pi\left(  g\right)  \right)  \in\wedge^{2}\mathfrak{g}$. Since
\[
L_{gh}^{-1}\left[  \pi\left(  gh\right)  -R_{h}\pi\left(  g\right)  \right]
=\tilde{\pi}\left(  gh\right)  -L_{h}^{-1}L_{g}^{-1}R_{h}\pi\left(  g\right)
=\tilde{\pi}\left(  gh\right)  -L_{h}^{-1}R_{h}L_{g}^{-1}\pi\left(  g\right)
\]%
\[
=\tilde{\pi}\left(  gh\right)  -L_{h}^{-1}R_{h}\tilde{\pi}\left(  g\right)
=\tilde{\pi}\left(  gh\right)  -\operatorname{Ad}_{h}^{-1}\left(  \tilde{\pi
}\left(  g\right)  \right)  ,
\]
a subgroup $H$ is $\pi$-coisotropic if and only if
\[
\tilde{\pi}\left(  gh\right)  -\operatorname{Ad}_{h}^{-1}\left(  \tilde{\pi
}\left(  g\right)  \right)  \in\mathfrak{h}\wedge\mathfrak{g}%
\]
for all $g\in G$ and $h\in H$. Thus $\operatorname{Ad}_{\sigma}H$ is $\pi
$-coisotropic if and only if
\[
\text{(*)}\;\;\tilde{\pi}\left(  g\operatorname{Ad}_{\sigma}\left(  h\right)
\right)  -\operatorname{Ad}_{\operatorname{Ad}_{\sigma}\left(  h\right)
}^{-1}\left(  \tilde{\pi}\left(  g\right)  \right)  \in\operatorname{Ad}%
_{\sigma}\left(  \mathfrak{h}\right)  \wedge\mathfrak{g}%
\]
for all $g\in G$ and $h\in H$. Similarly, $H$ is $\pi_{\sigma}$-coisotropic if
and only if
\[
\text{(**)}\;\;\widetilde{  \pi_{\sigma}}\left(  gh\right)
-\operatorname{Ad}_{h}^{-1}\left(  \widetilde{  \pi_{\sigma}}
\left(  g\right)  \right)  \in\mathfrak{h}\wedge\mathfrak{g}%
\]
for all $g\in G$ and $h\in H$. Note that
\[
 \widetilde{ \pi_{\sigma}}\left(  g\right)  =L_{g}%
^{-1}\left(  \pi_{\sigma}\left(  g\right)  \right)  =L_{g}^{-1}\left(
R_{\sigma}\pi\left(  g\sigma^{-1}\right)  \right)  =R_{\sigma}L_{g}^{-1}%
\pi\left(  g\sigma^{-1}\right)
\]%
\[
=R_{\sigma}L_{\sigma}^{-1}L_{\sigma}L_{g}^{-1}\pi\left(  g\sigma^{-1}\right)
=R_{\sigma}L_{\sigma}^{-1}L_{g\sigma^{-1}}^{-1}\pi\left(  g\sigma^{-1}\right)
=R_{\sigma}L_{\sigma}^{-1}\tilde{\pi}\left(  g\sigma^{-1}\right)
\]
\[
=\operatorname{Ad}_{\sigma}^{-1}\left(  \tilde{\pi}\left(  g\sigma
^{-1}\right)  \right)
\]
for any $g\in G$. So
\[
\widetilde{  \pi_{\sigma}}\left(  gh\right)
-\operatorname{Ad}_{h}^{-1}\left(  \widetilde{ \pi_{\sigma}}
\left(  g\right)  \right)  =\operatorname{Ad}_{\sigma}^{-1}\left(  \tilde
{\pi}\left(  gh\sigma^{-1}\right)  \right)  -\operatorname{Ad}_{h}^{-1}\left(
\operatorname{Ad}_{\sigma}^{-1}\left(  \tilde{\pi}\left(  g\sigma^{-1}\right)
\right)  \right)
\]%
\[
=\operatorname{Ad}_{\sigma}^{-1}\left(  \tilde{\pi}\left(  g\sigma
^{-1}\operatorname{Ad}_{\sigma}\left(  h\right)  \right)  \right)
-\operatorname{Ad}_{\sigma}^{-1}\operatorname{Ad}_{\sigma}\operatorname{Ad}%
_{h}^{-1}\operatorname{Ad}_{\sigma}^{-1}\left(  \tilde{\pi}\left(
g\sigma^{-1}\right)  \right)
\]%
\[
=\operatorname{Ad}_{\sigma}^{-1}\left[  \left(  \tilde{\pi}\left(
g\sigma^{-1}\operatorname{Ad}_{\sigma}\left(  h\right)  \right)  \right)
-\operatorname{Ad}_{\sigma h\sigma^{-1}}^{-1}\left(  \tilde{\pi}\left(
g\sigma^{-1}\right)  \right)  \right]
\]
and hence the condition (**) is equivalent to
\[
\tilde{\pi}\left(  g\sigma^{-1}\operatorname{Ad}_{\sigma}\left(  h\right)
\right)  -\operatorname{Ad}_{\operatorname{Ad}_{\sigma}h}^{-1}\left(
\tilde{\pi}\left(  g\sigma^{-1}\right)  \right)  \in\operatorname{Ad}_{\sigma
}\left(  \mathfrak{h}\wedge\mathfrak{g}\right)  =\operatorname{Ad}_{\sigma}\left(
\mathfrak{h}\right)  \wedge\mathfrak{g}%
\]
for all $g\in G$ and $h\in H$, or equivalently, the condition (*).
\fidi
\begin{proposition} Let $G$ be a Poisson--Lie group, $H$ a closed subgroup such
that its conjugate $Ad_{\sigma}H:=\operatorname{Ad}_{\sigma}H=\sigma
H\sigma^{-1}$ is coisotropic where $\sigma\in G$. Let $\pi_{\sigma}$ be the
affine Poisson structure on $G$ given by $\pi_{\sigma}(g):=R_{\sigma}%
\pi(g\sigma^{-1})$. Let $p:G\rightarrow G/H$ and $p_{\sigma}:G\rightarrow
G/H_{\sigma}$ be the natural projections. Then the Poisson manifolds
$(G/H,p_{\ast}\pi_{\sigma})$ and $(G/H_{\sigma},\left(  p_{\sigma}\right)
_{\ast}\pi)$ are Poisson diffeomorphic.
\end{proposition}\label{Poisson-diffeo}
\proof There is a natural diffeomorphism between $G/H$ and
$G/H_{\sigma}$ given by $\iota:\left[  g\right]  _{H}\mapsto\left[
g\sigma^{-1}\right]  _{H_{\sigma}}$ which satisfies
\[
p_{\sigma}=\iota\circ p\circ R_{\sigma},
\]
where $\left[  g\right]  _{H}:=gH\in G/H$. We claim that $\iota$ is a Poisson
map, i.e. $\iota_{\ast}\left(  p_{\ast}\pi_{\sigma}\right)  =\left(
p_{\sigma}\right)  _{\ast}\pi$. Indeed %

\[
\iota_{\ast}\left(  p_{\ast}\pi_{\sigma}\right)  \left(  \left[  g\sigma
^{-1}\right]  _{H_{\sigma}}\right)  =\left(  D\iota\right)  |_{\left[
g\right]  _{H}}\left(  (p_{\ast}\pi_{\sigma})\left(  \left[  g\right]
_{H}\right)  \right)  =\left(  D\iota\right)  |_{\left[  g\right]  _{H}%
}\left(  \left(  Dp\right)  |_{g}\left(  \pi_{\sigma}(g)\right)  \right)
\]%
\[
=\left(  D\iota\right)  |_{\left[  g\right]  _{H}}\left(  \left(  Dp\right)
|_{g}\left(  R_{\sigma}\pi(g\sigma^{-1})\right)  \right)  =D(\iota\circ p\circ
R_{\sigma})|_{g\sigma^{-1}}\left(  \pi(g\sigma^{-1})\right)
\]%

\[
=D(p_{\sigma})|_{g\sigma^{-1}}\left(  \pi(g\sigma^{-1})\right)  =\left(
(p_{\sigma})_{\ast}\pi\right)  (\left[  g\sigma^{-1}\right]  _{H_{\sigma}})
\]
for any $\left[  g\right]  _{H}\in G/H$ and the claim follows.
\fidi

\subsection{Foliation point of view}

It is somewhat
unexpected that the $\pi$-coisotropy of a conjugate subgroup
$\operatorname{Ad}_{\sigma}H$ is not related to the $\operatorname{Ad}%
_{\sigma}\left(  \pi\right)  $-coisotropy of the subgroup $H$ but related to
the $R_{\sigma}\pi$-coisotropy of $H$. In this section, we use a foliation
viewpoint to give a more conceptual explanation of this phenomenon.
We call a foliation $\mathcal{F}$ on a
manifold $M$ regular if the leaf (i.e. the quotient) space
$M/\mathcal{F}$ inherits a well-defined manifold structure from
$M$.

Let $\mathcal{F}$ be a regular foliation on a manifold $M $
and $\rho\in\wedge^{k}TM$ be a tensor field on $M$. We call
$\mathcal{F} $ \emph{$\rho $-coisotropic} if for any element $\left[
\eta\right]  $ of the holonomy groupoid $\mathfrak{G}$ that goes from
$s\in M$ to $t\in M$ (and hence $s,t$ belong to the same leaf $L$
of $\mathcal{F}$), there is a (leaf-preserving) local
diffeomorphism $\eta$, implementing $\left[  \eta\right]  $, from
a neighborhood of $s$ to a neighborhood of $t$ with $\eta\left(
s\right)  =t$, such that
\[ \left(  D\eta\right)
_{s}\left(  \rho\left(  s\right)  \right)  -\rho\left( t\right)
\in T_{t}L\wedge\left(  \wedge^{k-1}T_{t}M\right)  ,
\]
and hence
$\rho$ projects to a well-defined tensor field $\left[ \rho\right]
=\rho/\mathcal{F}$ on $M/\mathcal{F}$. Note that the differential
$D\eta$ of $\eta$ is a local vector bundle map from
$T\mathcal{F}$ to $T\mathcal{F}$, where
$T\mathcal{F}=\cup_{L\in\mathcal{F}}TL$.

Fix a tensor field $\pi$ on $G$. We consider the category
$\mathcal{C}$ of $G$-manifolds $M$ endowed with $\pi$-covariant
tensor field $\rho$ and a regular
$\rho$-coisotropic foliation $\mathcal{F}$ that is invariant under the $G$-action on $M$. A morphism between two
objects $\left(  M,\rho,\mathcal{F}\right)  $ and $\left(  \tilde{M}%
,\tilde{\rho},\mathcal{\tilde{F}}\right)  $ is a smooth
$G$-equivariant map $\phi:M\rightarrow\tilde{M}$, i.e. $\phi\left(
gx\right)  =g\phi\left( x\right)  $ for all $\left(  g,x\right)
\in G\times M$, that induces a well-defined smooth map $\left[
\phi\right]  :M/\mathcal{F}\rightarrow
\tilde{M}/\mathcal{\tilde{F}}$ and sends the tensor field $\rho$
on $M$ to $\tilde{\rho}$ on $M $, i.e. $\left(  D\phi\right)
\left(  \rho\left( x\right)  \right)  =\tilde{\rho}\left(
\phi\left(  x\right)  \right)  $ for all $x\in M$.

It is easily recognized that the map $\left[  \phi\right]  $ induced by such a
morphism $\phi$ is automatically $G$-equivariant and sends $\left[
\rho\right]  $ to $\left[  \tilde{\rho}\right]  $.
It is natural to see that for a given object $\left(  M,\rho,\mathcal{F}%
\right)  $ of $\mathcal{C}$, any diffeomorphism $\phi:M\rightarrow
M$ produces an object $\left(
\tilde{M},\tilde{\rho},\mathcal{\tilde{F}}\right)  $ of
$\mathcal{C}$ with $\tilde{M}=M$,
$\mathcal{\tilde{F}}=\phi_{\ast}\mathcal{F}$ whose leaves are
exactly the images of leaves of $\mathcal{F}$ under $\phi$, and
$\tilde{\rho}=\left(  D\phi\right)  \left(  \rho\right)  $, such
that $\phi$ becomes an invertible morphism from $\left(
M,\rho,\mathcal{F}\right)  $ to $\left(
\tilde{M}=M,\tilde{\rho},\mathcal{\tilde{F}}\right)  $. In
particular, $\tilde{\rho}=\left(  D\phi\right)  \left(
\rho\right)  $ is $\pi$-covariant just like $\rho$, and
$\mathcal{\tilde{F}}=\phi_{\ast }\mathcal{F}$ is
$\tilde{\rho}$-coisotropic.

For each connected closed subgroup $H$
of $G$, the $G$-manifold $G$ has a regular foliation
$\mathcal{F}_{H}$ with the right cosets $gH$, $g\in G$, as leaves,
such that each holonomy groupoid element $\left[  \eta\right]  $
is implemented by a right translation $R_{h}$ with $h\in H$, which
implies that for any tensor field $\rho$ on $G$, $\mathcal{F}_{H}$
is $\rho$-coisotropic if and only if the subgroup $H$ is
$\rho$-coisotropic. Note that the diffeomorphism
$R_{\sigma}:G\rightarrow G$ with $\sigma\in G$ maps the
foliation $\mathcal{F}_{H}$ determined by $H$ to the foliation $\mathcal{F}%
_{\operatorname{Ad}_{\sigma^{-1}}H}$ determined by $\operatorname{Ad}%
_{\sigma^{-1}}H$, because it sends the leaf $gH$ of
$\mathcal{F}_{H}$ to the leaf \[ \left(  gH\right)
\sigma=g\sigma\left(  \sigma^{-1}H\sigma\right)
=g\sigma\operatorname{Ad}_{\sigma^{-1}}H \] of
$\mathcal{F}_{\operatorname{Ad}_{\sigma^{-1}}H}$ for all $g\in G$.
Thus $R_{\sigma}$ determines an invertible morphism from $\left(
G,\rho
,\mathcal{F}_{H}\right)  $ to $\left(  G,R_{\sigma}\rho,\mathcal{F}%
_{\operatorname{Ad}_{\sigma^{-1}}H}\right)  $ for any tensor field
$\rho$ on $G$ that makes the subgroup $H$ $\rho$-coisotropic.
(This means that under the diffeomorphism $R_{\sigma}$, the tensor
field $\rho$ corresponds to
$R_{\sigma}\rho$ while the subgroup $H$ corresponds to $\operatorname{Ad}%
_{\sigma^{-1}}H$, not $R_{\sigma}H$ which is not a subgroup.)

In
particular, $\operatorname{Ad}_{\sigma^{-1}}H$ is
$R_{\sigma}\rho$-coisotropic when $H$ is $\rho$-coisotropic. Since
$R_{\sigma}$ is invertible, we have
$\operatorname{Ad}_{\sigma^{-1}}H$ is $R_{\sigma}\rho$-coisotropic
if and only if $H$ is $\rho$-coisotropic. Substituting $H$ by
$\operatorname{Ad}_{\sigma }H$, we can also say that $H$ is
$R_{\sigma}\rho$-coisotropic if and only if
$\operatorname{Ad}_{\sigma}H$ is $\rho$-coisotropic. Furthermore
from the above general discussion, it is also clear why the diffeomorphism
$R_{\sigma}$ induces a Poisson
diffeomorphism
$$\left(  G/H=G/\mathcal{F}_{H},\rho/\mathcal{F}%
_{H}\right)  \to \left(  G/\operatorname{Ad}_{\sigma^{-1}}H=G/\mathcal{F}%
_{\operatorname{Ad}_{\sigma^{-1}}H},\rho/\mathcal{F}_{\operatorname{Ad}%
_{\sigma^{-1}}H}\right)\, .  $$

\section{Poisson Grassmannians}

\subsection{Coisotropic subgroups in standard $SU(n)$}

Let us now restrict ourselves to the group $SU(n)$
and fix the embedding of $S(U(n-m)\times U(m))$ in $SU(n)$ given by:
$$
 (A,B)\hookrightarrow\begin{pmatrix}
   A & \underline{0} \\
   \underline{0} & B \
 \end{pmatrix}\, .
$$
Recall that the standard Poisson--Lie tensor on $SU(n)$ is defined, up to
a constant factor by the Poisson $2$--tensor
$$
\pi(g)=L_g r-R_g r
$$
where $r\in{\mathfrak g}\wedge{\mathfrak g}$, $\mathfrak g=\mathfrak{su}(n)$ in the following,
is the $r$-matrix given by
$$
r=\sum_{1\le i<j\le n}X^+_{ij}\wedge X^-_{ij}\, .
$$
Here we are considering the Cartan decomposition of $\mathfrak g$ determined by the
subalgebra of diagonal matrices and denote by $X^\pm_{ij}$ the
corresponding root vectors
$$
X^+_{ij}=\imath(e_{ij}+e_{ji})\qquad X^-_{ij}=e_{ij}-e_{ji}
$$
with $e_{ij}$ denoting a standard matrix unit.

It is then easily seen that $S(U(m)\times U(n-m))$ is a Poisson--Lie subgroup
in $SU(n)$. We will denote its Lie algebra by
$\mathfrak s(\mathfrak{u}(n-m)\times\mathfrak{u}(m))$.
\begin{proposition}
Let $m\leq\lbrack\frac{n}{2}]$ with $4\leq n\in\mathbb{N}$,
$c\in\lbrack0,1]$, and let
\[
\sigma(c,m)=\sqrt{c}\sum_{i=1}^{m}e_{ii}+\sum_{i=m+1}^{n-m}\! e_{ii}+\sqrt{c}\sum
_{i=n-m+1}^{n}e_{ii}+\sqrt{1-c}\sum_{i=1}^{m}(e_{n+1-i,i}-e_{i,n+1-i})\,.
\]
Then subgroup $Ad_{\sigma(c,m)}(U(n-m)\times U(m))$ is coisotropic in $U(n)$.
Analogously $Ad_{\sigma(c,m)}S(U(n-m)\times U(m))$ is coisotropic in the
standard $SU(n)$.
\end{proposition}
\noindent Proof
 Let $\sigma=\sigma(c,m)$ throughout the proof.
Then
\[
\sigma^{-1}=\sqrt{c}(\sum_{i=1}^{m}e_{i,i}+e_{n+1-i,n+1-i})+\sum_{k=m+1}%
^{n-m}\! e_{k,k}-\sqrt{1-c}(\sum_{i=1}^{m}e_{n+1-i,i}-e_{i,n+1-i})\,.
\]

As in the proof of Theorem 3 of \cite{Sh4}, it
suffices to show that
\[
\left(  Ad_{\sigma^{-1}}r\right)  -(2c-1)r\in\left(  \mathfrak{u}\left(
n-m\right)  \times\mathfrak{u}\left(  m\right)  \right)  \wedge\mathfrak{u}\left(
n\right)  .
\]

Let $A=\sqrt{c(1-c)}$. First of all we remark that the following relations
hold true:
\begin{align*}
\sigma^{-1}e_{i,j}\sigma & =ce_{i,j}+(1-c)e_{n+1-i,n+1-j}-A(e_{i,n+1-j}%
+e_{n+1-i,j})\\
\sigma^{-1}e_{n+1-j,n+1-i}\sigma & =ce_{n+1-j,n+1-i}+(1-c)e_{j,i}%
+A(e_{n+1-j,i}+e_{j,n+1-i})\\
\sigma^{-1}e_{i,n+1-j}\sigma & =ce_{i,n+1-j}-(1-c)e_{n+1-i,j}+A(e_{i,j}%
-e_{n+1-i,n+1-j})\\
\sigma^{-1}e_{n+1-i,j}\sigma & =ce_{n+1-i,j}-(1-c)e_{i,n+1-j}+A(e_{i,j}%
-e_{n+1-i,n+1-j})\\
&
\end{align*}
for every $1\leq i,j\leq m$. Furthermore:
\begin{align*}
\sigma^{-1}e_{i,m+p}\sigma & =\sqrt{c}e_{i,m+p}-\sqrt{1-c}e_{n+1-i,m+p}\\
\sigma^{-1}e_{n+1-i,m+p}\sigma & =\sqrt{c}e_{n+1-i,m+p}+\sqrt{1-c}e_{i,m+p}\\
\sigma^{-1}e_{m+p,i}\sigma & =\sqrt{c}e_{m+p,i}-\sqrt{1-c}e_{m+p,n+1-i}\\
\sigma^{-1}e_{m+p,n+1-i}\sigma & =\sqrt{c}e_{m+p,n+1-i}+\sqrt{1-c}e_{m+p,i}\\
&
\end{align*}
for every $1\leq i\leq m$ and $1\leq p\leq n-2m$. Lastly:%

\[
\sigma^{-1}e_{i,j}\sigma=e_{i,j}%
\]
when $m+1\leq i,j\leq n-m$. From these equalities one can compute:
\[
\sigma^{-1}X_{i,n+1-j}^{\pm}\sigma=cX_{i,n+1-j}^{\pm}\mp(1-c)X_{j,n+1-i}^{\pm
}+A(X_{i,j}^{\pm}\mp X_{n+1-j,n+1-i}^{\pm})
\]
for all $1\leq i\ne j\leq m$, and

\begin{align*}
\sigma^{-1}X_{i,j}^{\pm}\sigma & =cX_{i,j}^{\pm}-AX_{i,n+1-j}^{\pm}\mp
AX_{j,n+1-i}^{\pm}\pm(1-c)X_{n+1-j,n+1-i}^{\pm}\\
\sigma^{-1}X_{n+1-j,n+1-i}^{\pm}\sigma & =cX_{n+1-j,n+1-i}^{\pm}%
\pm(1-c)X_{i,j}^{\pm}\pm AX_{i,n+1-j}^{\pm}+AX_{j,n+1-i}^{\pm}\\
\sigma^{-1}X_{i,n+1-i}^{-}\sigma & =X_{i,n+1-i}^{-}\\
\sigma^{-1}X_{i,n+1-i}^{+}\sigma & =(2c-1)X_{i,n+1-i}^{+}+2AK_{i}%
\end{align*}
for every $1\leq i<j\leq m$, where $K_{i}=\imath(e_{i,i}%
-e_{n+1-i,n+1-i})$. Furthermore:%

\begin{align*}
\sigma^{-1}X_{i,m+p}^{\pm}\sigma & =\sqrt{c}X_{i,m+p}^{\pm}\mp\sqrt
{1-c}X_{m+p,n+1-i}^{\pm}\\
\sigma^{-1}X_{m+p,n+1-i}^{\pm}\sigma & =\sqrt{c}X_{m+p,n+1-i}^{\pm}\pm
\sqrt{1-c}X_{i,m+p}^{\pm}\\
&
\end{align*}
for every $1\leq i\leq m$ and $1\leq p\leq n-2m$. Lastly:%

\[
\sigma^{-1} X^{\pm}_{m+p,m+q}\sigma=X^{\pm}_{m+p,m+q}
\]
when $1\le p<q\le n-2m$.

Let's now move to $Ad_{\sigma^{-1}}r$ which we divide into three separate
pieces:
\[
Ad_{\sigma^{-1}}r=Ad_{\sigma^{-1}}(\Phi+\Theta+\Omega)
\]
where:%

\begin{align*}
\Phi & =\sum_{m+1\leq i,j\leq n-m}X_{i,j}^{+}\wedge X_{i,j}^{-}\\
\Theta & =\sum_{p=1}^{n-2m}\sum_{i=1}^{m}  X_{i,m+p}^{+}\wedge
X_{i,m+p}^{-}+X_{m+p,n+1-i}^{+}\wedge X_{m+p,n+1-i}^{-}+\sum
_{i=1}^{m}X_{i,n+1-i}^{+}\wedge X_{i,n+1-i}^{-}\\
\Omega & =\!\sum_{1\leq i<j\leq m}X_{i,j}^{+}\wedge X_{i,j}^{-}+\!\sum_{1\leq
i\neq j\leq m}\!\! X_{i,n+1-j}^{+}\wedge X_{i,n+1-j}^{-}+\!\!\sum_{1\leq i<j\leq
m}\!\! X_{n+1-j,n+1-i}^{+}\wedge X_{n+1-j,n+1-i}^{-}%
\end{align*}

By a straightforward computation, we get
\begin{align*}
Ad_{\sigma^{-1}}(\Phi)  & =\Phi
\in(2c-1)\Phi+\left[  \left(  \mathfrak{u}\left(  n-m\right)  \times
\mathfrak{u}\left(  m\right)  \right)  \wedge\mathfrak{u}\left(  n\right)  \right]  \\
Ad_{\sigma^{-1}}(\Theta)  & =(2c-1)\Theta+2A\sum[K_{i}\wedge X_{i,n+1-i}%
^{-}+X_{i,m+p}^{+}\wedge X_{m+p,n+1-i}^{-}-X_{m+p,n+1-i}^{+}\wedge
X_{i,m+p}^{-}]\\
& \in(2c-1)\Theta+\left[  \left(  \mathfrak{u}\left(  n-m\right)  \times
\mathfrak{u}\left(  m\right)  \right)  \wedge\mathfrak{u}\left(  n\right)  \right]
\end{align*}
since $K_{i},X_{i,m+p}^{\pm},\Phi\in\mathfrak{u}\left(  n-m\right)
\times\mathfrak{u}\left(  m\right)  $. The computation of $Ad_{\sigma
^{-1}}(\Omega)$ is much more tedious. It involves $Ad_{\sigma^{-1}}(\Omega)=$
\begin{align*}
 \sum_{1\leq i<j\leq m}\left(  cX_{i,j}^{+}-AX_{i,n+1-j}^{+}-AX_{j,n+1-i}
^{+}+(1-c)X_{n+1-j,n+1-i}^{+}\right) \\
 \wedge\left(  cX_{i,j}^{-}%
-AX_{i,n+1-j}^{-}+AX_{j,n+1-i}^{-}-(1-c)X_{n+1-j,n+1-i}^{-}\right)  \\
 +\sum_{1\leq i\neq j\leq m}\left(  AX_{i,j}^{+}+cX_{i,n+1-j}^{+}%
-(1-c)X_{j,n+1-i}^{+}-AX_{n+1-j,n+1-i}^{+}\right)\\
\wedge\left(  AX_{i,j}%
^{-}+cX_{i,n+1-j}^{-}+(1-c)X_{j,n+1-i}^{-}+AX_{n+1-j,n+1-i}^{-}\right)  \\
 +\sum_{1\leq i<j\leq m}\left(  (1-c)X_{i,j}^{+}+AX_{i,n+1-j}^{+}%
+AX_{j,n+1-i}^{+}+cX_{n+1-j,n+1-i}^{+}\right) \\
\qquad \wedge\left(  -(1-c)X_{i,j}%
^{-}-AX_{i,n+1-j}^{-}+AX_{j,n+1-i}^{-}+cX_{n+1-j,n+1-i}^{-}\right)  .
\end{align*}
The sum of all the wedge products of a ``$^{+}$--term'' on the
left of $\wedge$ and the corresponding ``$^{-}$-term'' on
the right of $\wedge$ is $(2c-1)\Omega$. All the remaining
wedge products of a term on the left of $\wedge$ and a term on the
right of $\wedge$ are in $\left(  \mathfrak{u}\left(  n-m\right)
\times\mathfrak{u}\left(  m\right)  \right)  \wedge u\left(  n\right)  $,
except for those involving the products $X_{i,n+1-j}^{+}\wedge X_{j,n+1-i}%
^{+}$ since multiples of $X_{i,n+1-j}^{+}$ and
$X_{j,n+1-i}^{+}$ are the only terms not in the Lie subalgebra
$u\left(  n-m\right)  \times u\left(  m\right)  $. It is easy to
check that the sum of all those wedge products involving $X_{i,n+1-j}%
^{+}\wedge X_{j,n+1-i}^{+}$ is $0$. So we get
\[
Ad_{\sigma^{-1}}(\Omega)\in(2c-1)\Omega+\left[  \left(  \mathfrak{u}\left(
n-m\right)  \times\mathfrak{u}\left(  m\right)  \right)  \wedge\mathfrak{u}\left(
n\right)  \right]
\]

Putting all together, we have that:
\[
Ad_{\sigma^{-1}}r=(2c-1)r+\left[  \left(  \mathfrak{u}\left(  n-m\right)
\times\mathfrak{u}\left(  m\right)  \right)  \wedge\mathfrak{u}\left(  n\right)
\right]
\]
as wanted.

\fidi

We will denote with $\tau_{\sigma_c}$ the projected Poisson 2--tensor
on the complex Grassmannian $G_{n}^{m}\mathbb{C}=SU(n)/S(U(m)\times
U(n-m))$.

\subsection{Covariance of tensor structures}

We plan now to describe a general argument which shows that the
Poisson pencil generated by the only (up to constant)
$SU(n)$--invariant Poisson structure on $G_{m}^{n}\mathbb{C}$
together with any $\tau_{\sigma_c}$ covers all of $SU(n)$--covariant
Poisson structures on the complex Grassmannians.

Let $M$ be a
$G$-manifold. Given two tensor fields $\pi$ and $\rho$ (of the
same kind) on $G$ and $M$ respectively, $\rho$ is called
$\pi$-covariant if (the differential $D\mu$ of) the action
\[
\mu:\left(  g,h\right)  \in G\times M\mapsto gh\in M
\]
sends the product tensor $\pi\times\rho$ on $G\times M$ to $\rho$ on $M$.
When $\pi,\rho$ are Poisson 2-tensors this means that
$\mu$ is a Poisson map (w.r.t. $\pi\times\rho$ and $\rho$), where the product tensor
$\pi\times\rho$ on $G\times M$ is defined by
\[
\left(\pi\times\rho\right)  \left(  g,h\right)  :=\pi\left(  g\right)
\oplus\rho\left(  h\right)  \in\wedge^{2}T_{g}G\oplus\wedge^{2}T_{h}%
M\subset\wedge^{2}T_{\left(  g,h\right)  }\left(  G\times M\right)\, .
\]
The general condition can be summarized as
\[
\left(  D\mu\right)\left(  \pi\oplus\rho\right)  =\rho\, .
\]
It is interesting to note that a tensor field $\rho$ on a $G$-manifold $M$ is $G$-invariant
if and only if $\rho$ is $0$-covariant for the vanishing tensor
field $0$ on $G$, i.e. the action operation
\[
\mu:\left(g,h\right)  \in G\times M\mapsto gh\in M
\]
sends the product tensor $0\times\rho$ on $G\times M$ to $\rho$ on $M$, because for
all $\left(  g_{0},h_{0}\right)  \in G\times M$, \[ \left(
D\mu\right)  _{\left(  g_{0},h_{0}\right)  }\left(  \left(
0\times\rho\right)  \left(  g_{0},h_{0}\right)  \right)  =\left(
D\mu\right) _{\left(  g_{0},h_{0}\right)  }\left(
0\oplus\rho\left(  h_{0}\right) \right) \] \[ =\left(
\frac{\partial\mu}{\partial g}\right)  _{\left(
g_{0},h_{0}\right) }\left(  0\right)  +\left(
\frac{\partial\mu}{\partial h}\right)  _{\left(
g_{0},h_{0}\right)  }\left(  \rho\left(  h_{0}\right)  \right)  =L_{g_{0}%
}\left(  \rho\left(  h_{0}\right)  \right) \] and hence \[ \left(
D\mu\right)  _{\left(  g_{0},h_{0}\right)  }\left(  \left(
0\times\rho\right)  \left(  g_{0},h_{0}\right)  \right)
=\rho\left( g_{0}h_{0}\right) \] if and only if \[
L_{g_{0}}\left(  \rho\left(  h_{0}\right)  \right)  =\rho\left(  g_{0}%
h_{0}\right)\,  . \]
Note that the multiplicativity of a tensor
field $\pi$ on $G$ is equivalent to the condition that $\pi$ is
$\pi$-covariant, i.e. the multiplication operation
\[ \mu:\left( g,h\right)  \in G\times G\mapsto gh\in G \]
sends the product tensor $\pi\times\pi$ on $G\times G$ to $\pi$ on $G$, because for
all $\left(  g_{0},h_{0}\right)  \in G\times G$, \[ \left(
D\mu\right)  _{\left(  g_{0},h_{0}\right)  }\left(  \left(  \pi
\times\pi\right)  \left(  g_{0},h_{0}\right)  \right)  =\left(
D\mu\right) _{\left(  g_{0},h_{0}\right)  }\left(  \pi\left(
g_{0}\right)  \oplus \pi\left(  h_{0}\right)  \right) \] \[
=\left(  \frac{\partial\mu}{\partial g}\right)  _{\left(
g_{0},h_{0}\right) }\left(  \pi\left(  g_{0}\right)  \right)
+\left(  \frac{\partial\mu }{\partial h}\right)  _{\left(
g_{0},h_{0}\right)  }\left(  \pi\left( h_{0}\right)  \right)
=R_{h_{0}}\left(  \pi\left(  g_{0}\right)  \right)
+L_{g_{0}}\left(  \pi\left(  h_{0}\right)  \right) \] and hence \[
\left(  D\mu\right)  _{\left(  g_{0},h_{0}\right)  }\left(  \left(
\pi
\times\pi\right)  \left(  g_{0},h_{0}\right)  \right)  =\pi\left(  g_{0}%
h_{0}\right) \] if and only if \[ R_{h_{0}}\left(  \pi\left(
g_{0}\right)  \right)  +L_{g_{0}}\left( \pi\left(  h_{0}\right)
\right)  =\pi\left(  g_{0}h_{0}\right)\,  . \]
Similarly, the affinity of a tensor field $\rho$ on $G$ is equivalent to the
condition that $\rho$ is $\pi$-covariant for the field $\pi=\rho_{l}%
:=\rho-\left(  \rho\left(  e\right)  \right)  ^{l}$ (which is
multiplicative when $\rho$ is indeed affine), i.e. the
multiplication operation \[ \mu:\left(  g,h\right)  \in G\times
G\mapsto gh\in G \] sends the product tensor $\pi\times\rho$ on
$G\times G$ to $\rho$ on $G$, because for all $\left(
g_{0},h_{0}\right)  \in G\times G$, \[ \left(  D\mu\right)
_{\left(  g_{0},h_{0}\right)  }\left(  \left(  \pi
\times\rho\right)  \left(  g_{0},h_{0}\right)  \right)  =\left(
D\mu\right) _{\left(  g_{0},h_{0}\right)  }\left(  \pi\left(
g_{0}\right)  \oplus \rho\left(  h_{0}\right)  \right) \] \[
=\left(  \frac{\partial\mu}{\partial g}\right)  _{\left(
g_{0},h_{0}\right) }\left(  \pi\left(  g_{0}\right)  \right)
+\left(  \frac{\partial\mu }{\partial h}\right)  _{\left(
g_{0},h_{0}\right)  }\left(  \rho\left( h_{0}\right)  \right)
=R_{h_{0}}\left(  \pi\left(  g_{0}\right)  \right)
+L_{g_{0}}\left(  \rho\left(  h_{0}\right)  \right) \] \[
=R_{h_{0}}\left(  \rho\left(  g_{0}\right)  \right)
-R_{h_{0}}\left( L_{g_{0}}\left(  \rho\left(  e\right)  \right)
\right)  +L_{g_{0}}\left( \rho\left(  h_{0}\right)  \right) \] and
hence \[ \left(  D\mu\right)  _{\left(  g_{0},h_{0}\right)
}\left(  \left(  \pi \times\rho\right)  \left(  g_{0},h_{0}\right)
\right)  =\rho\left( g_{0}h_{0}\right) \] if and only if \[
R_{h_{0}}\left(  \pi\left(  g_{0}\right)  \right)
+L_{g_{0}}\left( \rho\left(  h_{0}\right)  \right)
-R_{h_{0}}\left(  L_{g_{0}}\left( \rho\left(  e\right)  \right)
\right)  =\rho\left(  g_{0}h_{0}\right)\,  . \]
We give an interesting application of the above viewpoint. First we give a
proof of the following known general result, using the concept
discussed above.

\begin{proposition}. Let $\pi\in\wedge^{k}TG$ and $M$ be a
$G$-manifold. If $\rho\in\wedge^{k}TM$ is $\pi$-covariant and
$\tilde{\rho}\in\wedge^{k}TM$ is
$G$-invariant, then $\rho+\tilde{\rho}$ (or any tensor in $\rho+\mathbb{R}%
\tilde{\rho}$) is $\pi$-covariant and the Schouten bracket $\left[
\left[ \rho,\tilde{\rho}\right]  \right]  $ (or any tensor in
$\left[  \left[
\rho,\mathbb{R}\tilde{\rho}\right]  \right]  $) is $G$-invariant.%
\end{proposition}

\proof The given conditions can be summarized as \[ \left(
D\mu\right)  \left(  \pi\oplus\rho\right)  =\rho\;\text{and\
}\left( D\mu\right)  \left(  0\oplus\tilde{\rho}\right)
=\tilde{\rho} \] for the action map $\mu:G\times M\rightarrow M$.
Clearly we have \[ \left(  D\mu\right)  \left(  \pi\oplus\left(
\rho+\tilde{\rho}\right) \right)  =\left(  D\mu\right)  \left(
\left(  \pi\oplus\rho\right)  +\left( 0\oplus\tilde{\rho}\right)
\right)  =\rho+\tilde{\rho} \] which means that
$\rho+\tilde{\rho}$ is $\pi$-covariant. On the other hand, we
first note that the Schouten bracket \[
\left[  \left[  \kappa\oplus0,0\oplus\lambda\right]  \right]  =0\text{\ in\ }%
\wedge^{k}T\left(  G\times M\right) \] for any tensor
$\kappa\in\wedge^{k}TG$ and $\lambda\in\wedge^{k}TM$. Now since
the differential $D\mu$ preserves the Schouten bracket operation,
we also have \[ \left[  \left[  \rho,\tilde{\rho}\right]  \right]
=\left[  \left[  \left( D\mu\right)  \left(  \pi\oplus\rho\right)
,\left(  D\mu\right)  \left( 0\oplus\tilde{\rho}\right)  \right]
\right]  =\left(  D\mu\right)  \left( \left[  \left[
\pi\oplus\rho,0\oplus\tilde{\rho}\right]  \right]  \right) \] \[
=\left(  D\mu\right)  \left(  \left[  \left[  \pi,0\right]
\right] \oplus\left[  \left[  \rho,\tilde{\rho}\right]  \right]
\right)  =\left( D\mu\right)  \left(  0\oplus\left[  \left[
\rho,\tilde{\rho}\right]  \right] \right)  , \] which means that
$\left[  \left[  \rho,\tilde{\rho}\right]  \right]  $ is
$G$-invariant. \endproof

For Poisson tensors $\rho,\tilde{\rho}$ on $M$ (i.e. $\left[
\left[ \rho,\rho\right]  \right]  =0=\left[ \left[
\tilde{\rho},\tilde{\rho }\right]  \right]  $), the sum
$\rho+\tilde{\rho}$ is Poisson if and only if $\left[  \left[
\rho,\tilde{\rho}\right]  \right]  =0$. If $\rho\in\wedge ^{2}TM$
is a $\pi$-covariant Poisson tensor and
$\tilde{\rho}\in\wedge^{2}TM$ is $G$-invariant Poisson tensor,
then $\rho+\tilde{\rho}$ (or any tensor in
$\rho+\mathbb{R}\tilde{\rho}$) is a $\pi$-covariant Poisson tensor
if there is no non-trivial $G$-invariant 3-tensor on $M$.
For any compact symmetric space this last condition is equivalent to $H_{DR}^{3}\left(M\right)
=0$ which is verified, for example,  when $M=G_{m}^{n}
\mathbb{C}$. This proves that if $\rho$ is the $SU(n)$ invariant Poisson
tensor on $G_{m}^{n}\mathbb{C}$ then $\rho$ and $\tau_c$ are compatible ($\left[ \left[
\rho,\tau_c\right]\right]$) and therefore generates the Poisson pencil of
$SU(n)$--covariant Poisson tensors.

In particular, if $X\in\wedge^{2}\mathfrak{g}$ and $X^{l}$ is a
(of course $G$-invariant) Poisson 2-tensor on $G$, then
$\rho+X^{l}$ is an affine Poisson 2-tensor on $G$ (which is also
$\rho_{l}$-covariant and hence $\left( \rho+X^{l}\right)
_{l}=\rho_{l}$)  for any affine  Poisson 2-tensor $\rho$ on
$G$ (which is $\rho_{l}$-covariant for the multiplicative $\rho_{l}%
:=\rho-\left(  \rho\left(  e\right)  \right)  ^{l}$).

\subsection{Lagrangian subalgebras}

In \cite{Dr:phs} Drinfel'd showed how to relate Poisson homogeneous spaces
of a given Poisson--Lie group to orbits (under a natural action) of the group
itself on the algebraic
variety ${\cal L}$ of Lagrangian subalgebra of the double $D(\mathfrak
g)$. Such construction led Karolinsky (\cite{Ka:cph}) to a classification of
Poisson homogeneous spaces -- at least when $D(\mathfrak g)$ is complex semisimple -- in
terms of combinatorial data associated to the root system. Later on Evens and Lu
in \cite{ELu:vls} showed how to define a natural Poisson bivector on ${\cal
L}$ such that the Drinfel'd map is always an equivariant Poisson map. In this context
a quotient by a coisotropic subgroup corresponds to orbits in ${\cal L}$
containing at least one split subalgebra. In this paragraph we'll describe such Lagrangian
subalgebras for our specific family of covariant Poisson brackets on complex
Grassmannians, generalizing results in \cite{Sh4}.

\begin{lemma} Let $G$ be a Poisson--Lie group, $H$ a closed
connected subgroup, with Lie algebra $\mathfrak
h\subseteq\mathfrak g$. Let $\sigma\in G$ be such that ${\mathfrak
h}_\sigma=Ad_\sigma{\mathfrak h}$ is coisotropic in $\mathfrak g$
(i.e. $(d\pi)\big|_{{\mathfrak h}_\sigma}\subseteq {\mathfrak h}_\sigma\wedge \mathfrak g$).
Then the Lagrangian subalgebras corresponding to the Poisson
structure $\tau_\sigma$ on the homogeneous space $G/H$ over the
point $x_0:=eH$ is
 \begin{equation}\label{Lagrangian}
{\mathfrak h}+W={\mathfrak h}+\{ (x,\beta)\in {\mathfrak
g}\times{\mathfrak g}^*\, \big|\, \beta\in{\mathfrak h}^\bot,
\pi_{\sigma^{-1}}\lrcorner\beta =x\}
 \end{equation}
 \end{lemma}
\noindent Proof\par
 By construction (see \cite{Dr:phs}) the Lagrangian subalgebra over
 the point $\sigma\cdot x_0$ is split and equals ${\mathfrak l}={\mathfrak h}_\sigma\oplus{\mathfrak
 h}_\sigma^\bot$.
We will use $G$--equivariance of the correspondence between points and Lagrangian
subalgebras  Recall that  the action of $G$ on its double $D(\mathfrak g)\simeq \mathfrak g\oplus\mathfrak g^*$
 is given by
 $$
g\cdot (X,\alpha)=(Ad_g X+\tilde\pi(g)(Ad^*_g\alpha,-),Ad^*_g\alpha)
 $$
Therefore letting $g=\sigma^{-1}$ act on $(Ad_\sigma Y,Ad^*_\sigma \beta)$
 what we get is
 $$
\sigma\cdot {\mathfrak l}={\mathfrak h}\oplus\{ (x,\beta)\in
{\mathfrak g}\times{\mathfrak g}^*\, \big|\, \beta\in{\mathfrak
h}^\bot, \pi_{\sigma^{-1}}(e)\lrcorner\beta =x\}
 $$
Note that the Lagrangian subalgebra $\sigma\cdot\mathfrak l$ is
exactly the Lagrangian subalgebra complementary to $\mathfrak g$
associated with the affine Poisson bracket $\pi_\sigma$ (see
\cite{Lu}) \fidi

The Lagrangian subalgebras corresponding to
Poisson homogeneous complex Grassmannians, over the point $x_0=eH$ can be
computed either by solving
 $$
\pi_{\sigma^{-1}}(e)\lrcorner\beta = x
 $$
or remarking that ${\mathfrak h}_\sigma^\bot$ is generated, as a vector space, by the
following elements
\begin{eqnarray*}
&& \langle\mp x^{ij}_\pm+ x^{n+1-j,n+1-i}_\pm -\frac{2c-1}{\sqrt{c(1-c)}}x^{j,n+1-i}_\pm\qquad ,
1\le i<j\le m\rangle \\
&&\langle
x^{i,n+1-j}_\pm+\mp x^{j,n+1-i}_\pm\qquad , 1\le i<j\le m\rangle\\
&&
\langle\sqrt c x^{i,m+p}_\pm \pm\sqrt{1-c}
x^{m+p,n+1-i}_\pm\qquad , 1\le i\le m, 1\le p\le n-2m\rangle\\
&&\langle x^{i,n+1-i}_-,h_i+h_{n-i}+\frac{2c-1}{\sqrt{c(1-c)}}\sum _{j=n-i+1}^n
x^{n+1-j,j}_+\qquad , 1\le i\le m\rangle\end{eqnarray*}
where $x^{hk}_\pm$ are the dual elements of $X_{hk}^\pm$ and $h_l$
are the dual elements of the Cartan subalgebra standard basis
$H_l=\imath(e_{l,l}-e_{n,n})$.

Remark that ${\mathfrak h}_\sigma^\bot$ is a Lie subalgebra of
${\mathfrak g}^*$ and, as such, can be integrated to a
coisotropic subgroup $H^\bot$ of $G^*$. The Poisson homogeneous space $G^*/K^\bot$ is
called the \emph{complementary dual} of $G/H$ in \cite{CiGa} where it is shown that it fits into
a quantum duality scheme.

\section{Poisson embeddings}

\subsection{General embeddings}

\begin{lemma}\label{embed} Let $(G,\pi)$ be a Poisson--Lie group (with Lie cobracket $\delta$).
Let $K$ be a closed Poisson--Lie subgroup and let $H'$ be a closed coisotropic
subgroup in $G$. Then $H=K\cap H'$ is a coisotropic subgroup of
$K$ and the natural map $$
 \imath:K/H\to G/H';\, [k]_H\mapsto [k]_{H'}
$$ is a Poisson embedding with respect to the projected Poisson
structures. If $K\cup H'$ generates $G$ then $K/H$
is Poisson diffeomorphic to $G/H'$.\end{lemma}
 \noindent Proof\par
Coisotropy of $K\cap H'$  in $K$ follows from its infinitesimal characterization.
In fact, intersecting a subcoalgebra $\mathfrak k$ (i.e. $\delta(\mathfrak k)
\subseteq\mathfrak k\wedge\mathfrak k$) with a subcoideal $\mathfrak h'$
(i.e. $\delta(\mathfrak h')\subseteq\mathfrak h'\wedge\mathfrak g$) gives a subcoideal
of $\mathfrak k$.
The map $\imath$ is the unique map such that $p'\circ
i=\imath\circ p_H$, where $i:K\hookrightarrow G$ is the Poisson
embedding, $p':G\to G/H'$ and $p_H:K\to K/H$ are the
natural Poisson projections. It is then easily seen that $\imath$
is injective, Poisson and with injective differential.  This map is
also surjective if every $g\in G$ can be written as $g=kh'$, with
$h\in K$ and $h'\in H'$ so that the last statement follows as
well.
 \fidi

\vspace{0.5cm}

\noindent{\bf Examples:}
\begin{enumerate}
\item Let $K'\subseteq K\subseteq G$ be a chain of Poisson--Lie groups. Then the natural
map from $K/K'$ to $G/K'$ is a Poisson embedding. In this way, for
example, one can prove that standard Poisson spheres ${\mathbb
S}^{2k+1}$ (i.e. quotients $SU(k+1)/SU(k)$ w.r.t. the standard
Poisson $SU(k+1)$) are embedded in standard Poisson complex Stiefel
manifolds $V_{k}^n\mathbb C \simeq SU(n)/SU(k)$
 \item Let $H_\sigma$ be a 1--parameter family of coisotropic
subgroups containing a Poisson--Lie group $K=H_0$. Then we have a
Poisson embedding from $K/(K\cap H_\sigma)$ to $G/H_\sigma$. This
example will be frequently used in what follows.
 \item Let $G=SU(n)$ and let $K=SU(n-1)$ be the Poisson--Lie
subgroup of lower right corner matrices (i.e. the first row and
column are $(1,0,\ldots,0)$. Then let $H'=SU(n-1)$ be the
Poisson--Lie subgroup of upper left corner matrices. We have:
$H=SU(n-2)$ and $K/H\simeq {\mathbb S}^{2n-1}$ with the
standard Poisson structure, which is, then, naturally embedded in
$SU(n)/SU(n-1)\simeq{\mathbb S}^{2n+1}$. Taking $H_p=SU(n-p)$
as upper left corner matrices and repeating the argument we
find a chain of Poisson embeddings of spheres explaining the symplectic
foliation of the standard Poisson spheres.
\end{enumerate}

We will now give a description of Poisson embeddings for standard complex projective
spaces and complex Grassmannians and see how it relates with
the Bruhat-Poisson foliation. The same argument will then be generalized to
non standard complex Grassmannians (and projective spaces) in what follows.
\vspace{1cm}
\subsection{Complex projective spaces}

In this section the idea is to explain how the {\it subgroup
method} can be used to describe (part of) the symplectic foliation both for
standard and non standard complex projective spaces.
Let us recall that from the classification of Poisson--Lie subgroups of a given
standard compact Poisson--Lie group (see \cite{Stok}) one
can deduce that maximal Poisson--Lie subgroups in $SU(n)$ are the
diagonally embedded $S(U(k)\times U(n-k))$, $k=1,\ldots,n$.

Let us start with the \emph{standard} case. The complex projective
space $\mathbb{P}^n\mathbb{C}$ is identified with the quotient
$SU(n)/S(U(1)\times U(n-1))$ via the projection
\begin{eqnarray*}
& &p:SU(n)\rightarrow \mathbb{P}^{n-1}\mathbb{C}\\
& & A\mapsto [ A\cdot\,{}^t(0,\ldots,0,1)]=[A^{(n)}]
\end{eqnarray*}
where $A^{(i)}$ denotes the $i^{th}$--column of the matrix
$A$ and ${}^t(0,\ldots,0,1)$ the transposed column vector.
The corresponding standard Poisson structure has symplectic
foliation described by Schubert cells (see \cite{Stok} for more explicit description)
which is, in this case, described as a chain of embeddings
$$
\mathbb{P}^0\mathbb{C}\subseteq
\mathbb{P}^1\mathbb{C}\subseteq\ldots
\subseteq\mathbb{P}^{n-2}\mathbb{C}\subseteq\mathbb{P}^{n-1}\mathbb{C}
$$
each of which is given by equations $Z_1=\ldots Z_k=0$ in
homogeneous coordinates. It is then easily seen that the parabolic
subgroups corresponding to
$S=\{\alpha_1,\ldots,\widehat{\alpha_k},\ldots,\alpha_{n-1}\}$ intersects
$SU(n)$ in a Poisson-Lie subgroup
$$K_k=\left\{\left(\begin{array}{cc}
  A & 0  \\
  0 & B
\end{array}\right)\in S(U(k)\times U(n-k))\right\}$$
 having as image under the projection $p$ exactly $X_k=\mathbb{P}^{k-1}\mathbb{C}$.

\begin{theorem}
For any $k=1,\ldots,n-2$ we have
$K_k\cap K_{n-1}\simeq S(U(k)\times U(n-k-1)\times U(1))\simeq U(k)\times U(n-k-1)$. Furthermore
$K_k/K_k\cap K_{n-1}$ is Poisson diffeomorphic to the standard
Poisson $\mathbb{P}^{n-k-1}\mathbb{C}$ and projects onto $X_k$ via $p$.
\end{theorem}
\noindent\emph{Proof} \par
The statement about the intersection is easily verified. For the second statement
 consider the map
$$
\imath:\frac{SU(n-k)}{S(U(n-k-1)\times U(1))}\longrightarrow \frac{K_k}{K_k\cap
K_{n-1}};\qquad
\imath([B])=\left[\left(\begin{array}{cc}1&0\\
0&B\end{array}\right)\right]\, .
$$
This map is a Poisson diffeomorphism due to an application of
lemma \ref{embed} remarking that $\imath(S(U(n-k-1)\times U(1)))=K_k\cap SU(n-k)$, and
that the union $\left(K_{n-1}\cap K_k\right)\cup SU(n-k)$ generates  $K_k$.

\fidi

\noindent Note that $p(K_1)\supseteq\ldots p(K_{n-2})\supseteq
p(K_{n-1})=\{*\}$, i.e. all the embeddings granted by the
proposition are contained one into  another and overlap the
Schubert cell decomposition.

\vspace{0.3cm}

Let's move to the \emph{non standard case}. As we have seen in Proposition \ref{Poisson-diffeo} one can
consider it simply as obtained via a different projection, i.e.
identifying the complex projective space with a quotient of
$SU(n)$ as image of
\begin{eqnarray*}
& &p_\sigma:SU(n)\rightarrow \mathbb{P}^{n-1}\mathbb{C}\\
& &p_\sigma: A\mapsto [ A\cdot\,{}^t(\sqrt
c,0,\ldots,0,\sqrt{1-c})]=[\sqrt c A^{(1)}+\sqrt{1-c}\, A^{(n)}]\, .
\end{eqnarray*}
The stabilizer, in this case, is the subgroup $H_\sigma=Ad_{\sigma(c,1)} U(n-1)$.
Differently from the standard case, the Poisson--Lie subgroups $K_k$ have
images which are not contained one into another. In more detail
$p(K_k)$ consists of
$$
[\left(\begin{array}{cc}
  A & 0  \\
  0 & B
\end{array}\right)\cdot\,{}^t(\sqrt c,0,\ldots,0,\sqrt{1-c})]=[(\sqrt c a_{11},\ldots,\sqrt c a_{k1},
\sqrt{1-c}b_{k+1,n},\ldots,\sqrt{1-c}b_{n,n})]
$$
The images $X_k=p(K_k)§$ satisfy, then, the equation
$$
\|A^{(1)}\|^2-\|B^{(n)}\|^2=0
$$
which, in homogeneous coordinates, can be expressed as
$$
|Z_{11}|^2+\ldots +|Z_{k1}|^2-\frac{c}{(1-c)}(|Z_{k+1,n}|^2+\ldots
+|Z_{nn}|^2)=0
$$
These are exactly the same equations for the higher dimensional singular symplectic leaves as in
\cite{KRR}.

\begin{theorem}\label{nonstandardproj}
For any $k=1,\ldots,n-1$ we have
$$
K_k\cap H_\sigma=\left\{\left(\begin{array}{cccc}
a&0&0&0\\ 0&B_{11}&0&0\\ 0&0&B_{22}&0\\ 0&0&0&a
\end{array}\right)\,: a^2=\det(B_{11}B_{22})^{-1}, B_{11}\in U(k-1), B_{22}\in U(n-k-1)\right\}\, .
$$
Furthermore $K_k/K_k\cap H_{\sigma_c}$ is Poisson diffeomorphic to
$({\mathbb S}^{2k-1}\times{\mathbb S}^{2(n-k)-1})/{\mathbb T}$
(if $k=1$ to a standard Poisson odd sphere ${\mathbb S}^{2n-3}$).
\end{theorem}
\noindent\emph{Proof} Let us start with $k=1$ and
consider the embedding
\[ A\in U\left(  n-1\right)  \mapsto\left( \begin{array}
[c]{cc}%
\det\left(  A\right)  ^{-1} & 0\\ 0 & A \end{array} \right)  \in
H:=K_1\subset G:=SU\left(  n\right) \] of $U\left(  n-1\right)  $
\textbf{ONTO} the closed subgroup $H=K_1$ of $G=SU\left( n\right)  $.
Let
\[ \sigma_{c}=\left( \begin{array}
[c]{ccc}%
\sqrt{c} &0& \sqrt{1-c}\\ 0 & I_{n-2} & 0\\
-\sqrt{1-c} & 0 & \sqrt c%
\end{array} \right)  \in SU\left(  n\right) \]
with $c\in\left(0,1\right)  $. Since for any
\[
h:=\left( \begin{array}[c]{ccc}%
a & 0 & 0\\ 0 & B & C\\ 0 & D & b \end{array} \right)  \in H
\]
with $a,b\in\mathbb{C}$, the conjugate
\[
\sigma_{c}h\sigma_{c}^{-1}=\left( \begin{array}[c]{ccc}
\sqrt{c} &0& \sqrt{1-c}\\ 0 & I_{n-2} & 0\\
-\sqrt{1-c} & 0 & \sqrt c
\end{array} \right)
\left( \begin{array}
[c]{ccc}%
a & 0 & 0\\ 0 & B & C\\ 0 & D & b \end{array} \right)
\left(\begin{array}
[c]{ccc}%
\sqrt{c} &0& -\sqrt{1-c}\\ 0 & I_{n-2} & 0\\
\sqrt{1-c} & 0 & \sqrt c
\end{array} \right)
\]%
\[ =\left( \begin{array}
[c]{ccc}%
a\sqrt{c} & \fbox{$D\sqrt{1-c}$} & b\sqrt{1-c}\\ 0 &
B & C\\ -a\sqrt{1-c} & \sqrt c D & b\sqrt c \end{array} \right)  \left(
\begin{array}
[c]{ccc}%
\sqrt{c} &0& -\sqrt{1-c}\\ 0 & I_{n-2} & 0\\
\sqrt{1-c} & 0 & \sqrt c
\end{array} \right)
\]%
\[ =\left( \begin{array}
[c]{ccc}%
ac+b-bc & \fbox{$D\sqrt{1-c}$}&\left(  b-a\right)  \sqrt{c}\sqrt{1-c}
\\  \fbox{$C\sqrt{1-c}$} &
B & C\sqrt{c}\\ \left(  b-a\right)  \sqrt{c}\sqrt{1-c}& D\sqrt{c} &
a-ac+bc\end{array} \right)
\]
is in $H$ if and only if $C=0$, $D=0$, and $b=a$, in which case
\[
\sigma_{c}h\sigma_{c}^{-1}=\left(\begin{array}
[c]{ccc}%
a & 0 & 0\\ 0 & B & 0\\ 0 & 0 & a \end{array} \right)  =h
\]
with $B\in U\left(  n-2\right)  $ and $a^{2}=\det\left(  B\right)
^{-1}$. Thus
\[
K_{\sigma_{c}}:=H\cap H_{\sigma_{c}}=\left\{
\left( \begin{array}[c]{ccc}%
a & 0 & 0\\ 0 & B & 0\\ 0 & 0 & a \end{array} \right)  \in
SU\left(  n\right)  :B\in U\left(  n-2\right)  \right\}
\]%
\[
=\left\{  \left( \begin{array}[c]{ccc}%
a & 0 & 0\\ 0 & B & 0\\ 0 & 0 & a \end{array} \right)  :B\in
U\left(  n-2\right)  \;\text{and\ }a^{2}=\det\left(  B\right)
^{-1}\right\}
\]
is a double covering $\tilde{U}\left(  n-2\right)$ of $U\left(  n-2\right)$ where
$H_{\sigma_{c}}:=\operatorname{Ad}_{\sigma_{c}}H=\sigma_{c}H\sigma
_{c}^{-1}$, and
\[
\mathfrak{k}_{\sigma_{c}}=\left\{  \left(\begin{array}[c]{ccc}%
a & 0 & 0\\ 0 & B & 0\\ 0 & 0 & a \end{array}
\right)  :B\in\mathfrak{u}\left(  n-2\right)  \;\text{and\ }-2a=\operatorname{tr}%
\left(  B\right)  \right\}  \cong\mathfrak{u}\left(  n-2\right)  . \]
It is not immediately clear that
$H/K_{\sigma_{c}}\cong\mathbb{S}^{2n-3}$ since
$K_{\sigma_{c}}\cong\tilde{U}\left(  n-2\right)  \neq U\left(
n-2\right)  $ and furthermore under the following identification
of $H$ and $U\left(  n-1\right)  $, $K_{\sigma_{c}}$ is not
identified with the standard canonically embedded $U\left(
n-2\right)  $, namely, $\left\{  \left( \begin{array}
[c]{cc}%
1 & 0\\ 0 & B \end{array} \right)  :B\in U\left(  n-2\right)
\right\}  $.

Let us prove that
$H/K_{\sigma_{c}}\cong \mathbb{S}^{2n-3}$ \textbf{and} that $\pi$
on $H$ projects to the standard covariant Poisson structure on
$\mathbb{S}^{2n-3}$. Indeed since $K_{\sigma_{c}}$ is a
$\pi$-coisotropic subgroup of $H$ and the canonically embedded
\[
SU\left(  n-1\right)  \equiv H_{0}:=\left\{  \left( \begin{array}
[c]{cc}%
1 & 0\\ 0 & A \end{array} \right)  :A\in SU\left(  n-1\right)
\right\}
\]
in $H$ is a Poisson-Lie subgroup of $H$, we have a Poisson embedding
\[
\iota:H_{0}/\left(  H_{0}\cap
K_{\sigma_{c}}\right)  \rightarrow
H/K_{\sigma_{c}}%
\]
where the Poisson structures are projected from $\pi$. Note
that $\iota$ is surjective (and hence is a diffeomorphism) since
$H_{0}\cup K_{\sigma_{c}}$ generates the group $H$. Note also that
\[
H_{0}\cap K_{\sigma_{c}}=\left\{  \left( \begin{array}
[c]{ccc}%
1 & 0 & 0\\ 0 & B & 0\\ 0 & 0 & 1 \end{array} \right)  :B\in
SU\left(  n-2\right)  \right\}
\]
the canonically embedded $SU\left(  n-2\right)  $ in $SU\left(  n-1\right)  $ and hence
\[
H_{0}/\left(  H_{0}\cap K_{\sigma_{c}}\right)  =SU\left(
n-1\right) /SU\left(  n-2\right)  =\mathbb{S}^{2n-3}\, .
\]
This shows that $H/K_{\sigma_{c}}\cong H_{0}/\left(  H_{0}\cap K_{\sigma_{c}%
}\right)  =\mathbb{S}^{2n-3}$ the standard covariant Poisson
sphere. Recall that
$H_{\sigma_{c}}=\operatorname{Ad}_{\sigma_{c}}H$ is $\pi
$-coisotropic and $\pi$ projects to the non-standard covariant
Poisson structure on $\mathbb{C}P^{n-1}=G/H_{\sigma_{c}}$. So with
$H_{\sigma_c}$ being a Poisson-Lie subgroup of $\left(  G,\pi\right)
$, we have a Poisson embedding
\[
H/K_{\sigma_{c}}=\mathbb{S}^{2n-3}\rightarrow G/H_{\sigma_{c}}=\mathbb{C}%
P^{n-1}%
\]
of the standard Poisson $\mathbb{S}^{2n-3}$ into the non-standard Poisson $\mathbb{C}P^{n-1}$.
\par
Let now $k\ne 1$. We want to prove that
$$
\frac{K_k}{K_k\cap H_{\sigma_c}}\simeq \frac{{\mathbb S}^{2k-1}\times {\mathbb S}^{2(n-k)-1}}{\mathbb T}
$$
as Poisson manifold, clarifying which is the Poisson structure on the right. Repeating the same argument as in
the first part of the proof we easily see that
$K_k\cap H_{\sigma_c}$ consists of matrices
$$
\left(\begin{array}{cccc}a&0&0&0\\ 0&B_{11}&0&0\\ 0&0&B_{22}& 0\\ 0&0&0& a\end{array}\right)
$$
where $B_{11}\in U(k-1)$, $B_{22}\in U(n-k-1)$ and $a^2=det(B_{11}B_{22})^{-1}$.
Now since $K_k\cup J$ generates $U(k)\times U(n-k)$ we have
$$
\frac{K_k}{K_k\cap H_{\sigma_c}}\simeq \frac{U(k)\times U(n-k)}{J}
$$
where, as Poisson manifolds, $U(k)\times U(n-k)$ has the product Poisson structure (of standard Poisson $U(i)$'s)
and $J$ consists of matrices
$$
\left(\begin{array}{cccc}a&0&0&0\\ 0&B_{11}&0&0\\ 0&0&B_{22}& 0\\ 0&0&0& a\end{array}\right)
$$
with $a\in U(1)$, $B_{11}\in U(k-1)$, $B_{22}\in U(n-k-1)$ (hence $J$ is a Poisson--Lie subgroup
of $U(k)\times U(n-k)$).
We remark that
$$
\frac{U(k)\times U(n-k)}{1\times U(k-1)\times U(n-k-1)\times 1}\simeq {\mathbb S}^{2k-1}\times {\mathbb S}^{2(n-k)-1}
$$
with the product of standard Poisson structures on the right. It is just a quotient by a Poisson--Lie subgroup
of $U(k)\times U(n-k)$.

The canonical embedding
$$
1\times U(k-1)\times U(n-k-1)\times 1\subseteq J
$$
of Poisson--Lie groups induces a Poisson quotient map
$$
\frac{U(k)\times U(n-k)}{1\times U(k-1)\times U(n-k-1)\times 1}\twoheadrightarrow \frac{U(k)\times U(n-k)}J\, .
$$
Since the actions of the subgroups $1\times U(k-1)\times U(n-k-1)\times 1$ and $\mathbb T=\{a\oplus I_{k-1}\oplus
I_{n-k-1}\oplus a\, : a\in U(1)\}$ commute, $\mathbb T$ gives a well defined diagonal action on
$\frac{U(k)\times U(n-k)}{1\times U(k-1)\times U(n-k-1)\times 1}\simeq \mathbb S^{2k-1}\times\mathbb S^{2(n-k)-1}$
such that the quotient map onto its orbit space coincides with the above quotient map.

The symplectic foliation of the standard covariant Poisson $\mathbb S^{2k-1}$ consists of $\mathbb T$--families of
$\mathbb C^i$ for $0\le i\le k-1$ with the $\mathbb T$--action on $\mathbb S^{2k-1}$ taking a leaf $\mathbb C^i$ to a leaf
$\mathbb C^i$ in the same $\mathbb T$--family. So the symplectic foliation of $\mathbb S^{2k-1}\times \mathbb S^{2(n-k)-1}$
consists of $\mathbb T^2$--families of $\mathbb C^i\times \mathbb C^j$ for $0\le j\le n-k-1$ and hence the symplectic
foliation of $\mathbb S^{2k-1}\times \mathbb S^{2(n-k)-1}/\mathbb T$ consists of $\mathbb T$--families of
$\mathbb C^i\times\mathbb C^j$ for $0\le i\le k-1$ and $0\le j\le n-k-1$.
\fidi
\vspace{0.5cm}

\emph{\underline{Remarks}}
\begin{enumerate}
\item Note that $\dim X_k=\dim K_k-\dim (K_k\cap H_\sigma)=[k^2+(n-k)^2-1]-[(k-1)^2+(n-k-1)^2]=2n-3$
independently of $k$.
\item It is obvious that whenever $k\ne l$,
$X_k\cap X_l$ is a union of lower dimensional symplectic leaves.
Each such intersection is just the image under the Poisson embedding
of the Poisson--Lie subgroup $K_k\cap K_l$.
\item The embedding $i_\sigma$ is the same as the Poisson map of Theorem 5 in
\cite{Sh4}. To prove this statement consider that the map granted by proposition \ref{embed} can be
constructed as follows: take $(v_1,\ldots,v_{n-1})$ complex
coordinates on the sphere, take  $u'\in U(n-1)$ with last column
equal to $(v_1,\ldots, v_{n-1})$ and consider $1\oplus u'$ as the
matrix with first row and first column equal to $(1,0,\ldots,0)$.
Projecting this matrix with respect to $H_\sigma$ means projecting
with $p\circ R_\sigma$ so that a direct computation shows that
the Poisson map of proposition \ref{embed} is:
$$
 (v_1,\ldots, v_{n-1})\mapsto [\sqrt{1-c}, \sqrt c v_1,\ldots, \sqrt c v_{n-1}]
$$
(here $[.]$ stands for equivalence class in ${\mathbb P}^{n-1}$)
which is exactly the same map as in \cite{Sh4} (apart from composition with the obvious Poisson diffeomorphism
$c\to 1-c$). It is remarkable that the connected components of the
 complementary of the union of the images of such embeddings
 are exactly the Poisson leaves of higher rank. Furthermore lower
 dimensional leaves can also be described as intersections of a
 suitable number of such images (the intersection of Poisson
 submanifolds being again a Poisson manifold), so that one can, in
 fact, completely describe the symplectic foliation of the complex
 projective space.
\end{enumerate}

\subsection{Complex Grassmannians}

In this section we study the more general Grassmannian case.
Let us fix once and for all the complete flag in ${\mathbb C}^n$,
$V_i=\langle e_{n-i},\ldots,e_n\rangle$ and let us give
notations for the Schubert cell decomposition.
Let $(a_1,\ldots,a_k)$ be a $k$--tuple of integers such that
$0\le a_1\le\ldots\le a_k\le n-k$,
and denote with $[a_1,\ldots,a_k]$ the corresponding
Schubert cell, i.e. the
set of $k$--planes in ${\mathbb C}^n$:
$$
[a_1,\ldots,a_k]=\{X\in G_k^n\mathbb C\,\big|\, \dim(X\cap V_{a_i+i})\ge i\}
$$

Then $[a_1,\ldots,a_k]$ is a cell of complex dimension
$\sum_{i=1}^ka_i$. The relative position of cells is described by
the so called Bruhat order:
$$
(a_1,\ldots,a_k)\le (b_1,\ldots,
b_k)\Longleftrightarrow a_i\le b_i\quad\forall i=1,\ldots,k\, .
$$
This is a partial ordering on the $k$--tuples of integers such that $(a_1,\ldots a_k)
\le (b_1,\ldots,b_k)$ if and only if $[a_1,\ldots,a_k]\subseteq [b_1,\ldots, b_k]$.
 Notice that $[a_1,\ldots,a_{k-1}]\mapsto
[0,a_1,\ldots,a_{k-1}]$ describes an embedding of $G_{k-1}^n\mathbb C$ into
$G_k^n\mathbb C$.

Now we consider subgroups and their projections, starting with the \emph{standard}
case.

\begin{theorem}
For any $l=1,\ldots,n-1$, let $K_l=S(U(l)\times U(n-l))$ and let $G=SU(n)$. Then we have:
\begin{enumerate}
\item There is a Poisson diffeomorphism
$$
\frac{K_l}{K_l\cap K_k}=
  \begin{cases}
    G_{k-l}^{n-l}\mathbb C & \text{if}\quad l<k, \\
    \{e\} &\text{if} \quad l=k\\
    G_{k}^{l}\mathbb C & \text{otherwise}.
  \end{cases}\, ;
$$
\item The image $X_l:=p(K_l)$ of $K_l$ under the projection $p:G\to G/K_k$ is the submanifold
$$
X_l=
  \begin{cases}
 \left[  \underset{l}{\underbrace{0,...,0}},\underset{k-l}%
{\underbrace{n-k,...,n-k}}\right] & \text{if}\quad l<k, \\
 [0,\ldots,0]&\text{if}\quad l=k\\
    \left[  \underset{k}{\underbrace
{l-k,...,l-k}}\right]  & \text{otherwise}.
  \end{cases}\, .
$$
 Note that we have the following
inclusion relations: $X_1\supseteq\ldots\supseteq X_{k-1}$ and
$X_{k+1}\subseteq\ldots\subseteq X_{n-1}$.
\end{enumerate}
\end{theorem}
\noindent\emph{Proof}

 First we note that
\[
K_{l}\cap K_{k}=\left\{%
\begin{array}
[c]{ll}%
S(U(l)\times U(k-l)\times U(n-k)) & \text{if}\quad l<k,\\
K_{k} & \text{if}\quad l=k\\
S(U(k)\times U(l-k)\times U(n-l)) & \text{otherwise}.
\end{array}\right.
\]
Furthermore the union of the subgroups $1_{l}\times SU(n-l)$ and $K_{l}\cap
K_{k}$ generates $K_{l}$ with
\[
\left(  1_{l}\times SU(n-l)\right)  \cap K_{l}\cap K_{k}=1_{l}\times
S(U(k-l)\times U(n-k))
\]
if $l\le k$, while the union of the subgroups $SU(l)\times1_{n-l}$ and $K_{l}\cap
K_{k}$ generates $K_{l}$ with
\[
\left(  SU(l)\times1_{n-l}\right)  \cap K_{l}\cap K_{k}=S(U(k)\times
U(l-k))\times1_{n-l}%
\]
if $l>k$. So by lemma \ref{embed}, we get Poisson diffeomorphisms%
\begin{align*}
\frac{SU(n-l)}{S(U(k-l)\times U(n-k))}\rightarrow\frac{K_{l}}{K_{l}\cap K_{k}%
};\quad\lbrack B] &  \mapsto\left(
\begin{array}
[c]{cc}%
I_{l} & 0\\
0 & B
\end{array}
\right)  \qquad l\leq k\\
\frac{SU(l)}{S(U(k)\times U(l-k))}\rightarrow\frac{K_{l}}{K_{l}\cap K_{k}%
};\quad\lbrack B] &  \mapsto\left(
\begin{array}
[c]{cc}%
B & 0\\
0 & I_{n-l}%
\end{array}
\right)  \qquad l>k\,.
\end{align*}
The rest of the theorem comes from direct computations. \fidi

We remark that different from the case of complex projective spaces, Poisson
embeddings of lower dimensional homogeneous spaces do not cover the whole
symplectic foliation for the complex Grassmannians which coincides with the
Schubert cell decomposition.

Let us move to the \emph{non standard} situation. We are then
considering
$$
G_k^n\mathbb C\simeq SU(n)/Ad_{\sigma(c,k)} S(U(k)\times U(n-k))
$$
with the projected Poisson tensor $\tau_{\sigma_c}$. Let us
consider the family of maximal Poisson--Lie subgroups $S(U(l)\times U(n-l))$, $1\le l\le n-1$.
The problem is to describe, for every  $l$, the image of $S(U(l)\times U(n-l))$ in $G_k^n\mathbb C$
and the Poisson manifold
$$
S(U(l)\times U(n-l))/(S(U(l)\times
U(n-l))\cap Ad_{\sigma(c,k)} (S(U(k)\times U(n-k)))\, .
$$

Let  ${\mathbb J}_k$ denote the $k\times k$ anti--diagonal matrix
$$
{\mathbb J}_k=\left(\begin{array}{cccc}0&\ldots&0&1\\
0&\ldots&1&0\\ \vdots&\ddots&&\vdots\\ 1&0&\ldots&0\end{array}\right)=\sum_{i=1}^k e_{i,k-i+1}\, .
$$
In the following the subscript of $\mathbb J_k$ is often omitted since the size of $\mathbb J$ is varying and can be easily determined from its surrounding context. With this notation
$$
\sigma(c,k)=\left(\begin{array}{ccc}\sqrt c{\mathbb I}_k&0&-\sqrt{1-c}{\mathbb J}_k\\
0&{\mathbb I}_{n-2k}&0\\ \sqrt{1-c}{\mathbb J}_k&0&\sqrt c{\mathbb I}_k\end{array}\right)\, .
$$
\begin{lemma}\label{symmetry}
We have $K_l\cap Ad_{\sigma(c,k)} K_k=K_{n-l}\cap Ad_{\sigma(1-c,k)} K_k$.
\end{lemma}
\noindent\emph{Proof}
First of all  $Ad_{\sigma(c,k)} K_k$ consists of
matrices of the form
$$
\sigma(c,k)^{-1}
\left(\begin{array}{ccc}A&0&0\\ 0&B_{11}&B_{12}\\ 0&B_{21}&B_{22}\end{array}\right)
\sigma(c,k)
$$
$$
=\left(\begin{array}{ccc}cA+(1-c){\mathbb J}B_{22}{\mathbb J}&\sqrt{1-c}{\mathbb J}B_{21}
&\sqrt{c(1-c)}(-A{\mathbb J}+{\mathbb J}B_{22})\\
\sqrt{1-c}B_{12}{\mathbb J}&B_{11}&\sqrt c B_{12}\\
\sqrt{c(1-c)}(-{\mathbb J}A+B_{22}{\mathbb J})&\sqrt{c}B_{21}&
(1-c){\mathbb J}A{\mathbb J}+cB_{22}\end{array}\right)
$$

Now the main point is to remark that $K_{n-i}=Ad_{\mathbb J}K_i$
and that
$$
{\mathbb J}\left(\begin{array}{ccc}cA+(1-c){\mathbb J}B_{22}{\mathbb J}
&\sqrt{1-c}{\mathbb J}B_{21}
&\sqrt{c(1-c)}(-A{\mathbb J}+{\mathbb J}B_{22})\\
\sqrt{1-c}B_{12}{\mathbb J}&B_{11}&\sqrt c B_{12}\\
\sqrt{c(1-c)}(-{\mathbb J}A+B_{22}{\mathbb J})&\sqrt{c}B_{21}&
(1-c){\mathbb J}A{\mathbb J}+cB_{22}\end{array}\right)
{\mathbb J}
$$
$$
=\left(\begin{array}{ccc}(1-c)A+c{\mathbb J}B_{22}{\mathbb J}
&\sqrt{c}{\mathbb J}B_{21}{\mathbb J}
&\sqrt{c(1-c)}(-A{\mathbb J}+{\mathbb J}B_{22})\\
\sqrt{c}{\mathbb J}B_{12}{\mathbb J}&{\mathbb J}B_{11}{\mathbb J}&\sqrt{1-c}{\mathbb J}B_{12}\\
\sqrt{c(1-c)}(-{\mathbb J}A+B_{22}{\mathbb J})&\sqrt{1-c}B_{21}{\mathbb J}&
c{\mathbb J}A{\mathbb J}+(1-c)B_{22}\end{array}\right)
$$
$$
=Ad_{\sigma(1-c,k)}\left(\begin{array}{ccc}A&0&0\\ 0&{\mathbb J}B_{11}{\mathbb
J}& {\mathbb J}B_{12}\\ 0& B_{21} {\mathbb J} &
B_{22}\end{array}\right)
$$
From this the claim follows.
\fidi

\begin{theorem}
For any $l=1,\ldots,n-1$, let $K_l=S(U(l)\times U(n-l))$ and let
$X_{l,k}=K_l/(K_l\cap Ad_{\sigma(c,k)} K_k)$. Then we have:
\begin{enumerate}
\item If $l<k$ or $l>n-k$ then $X_{l,k}$ is Poisson diffeomorphic to
$$
SU(n-l)/\left(S(U(|k-l|)\times U(|n-k-l|))\times  1_l\right)
$$
with a non standard Poisson quotient structure. The image of $X_{l,k}$ is a
Poisson submanifold of $G_k^n\mathbb C$ of codimension $l^2$ if $l<k$
and of codimension $(n-l)^2$ if $l> n-k$.
\item If $l=k$ or $l=n-k$ then $X_{l,k}$ is Poisson diffeomorphic to the Stiefel manifold
$$
   V_{k}^{n-k}\mathbb C
$$
with the standard quotient structure; its image in $G_k^n\mathbb C$ is therefore
a Poisson submanifold of codimension $k^2$. Remark that when $2k=n$, $V_k^{k}\mathbb C\simeq U(k)$
(with the standard Poisson structure).
\item If $k<l<n-k$ then $X_{l,k}$ is Poisson diffeomorphic to the quotient space
$$
 \frac{V_k^l\mathbb C\times V_{k}^{n-l}\mathbb C}{U(k)}
$$
of the standard Poisson complex Stiefel manifold $V_k^l\mathbb C\times V_{k}^{n-l}\mathbb C$
by the diagonal action of $U(k)$.
 The image of the projection $X_{l,k}$ is a submanifold of $G_k^n\mathbb C$
of codimension
$k^2$.
\end{enumerate}
\end{theorem}
\noindent\emph{Proof}

 Let us start with the case $l\le k$. From the formula for $Ad_{\sigma(c,k)}K_k$ described in the proof of the
preceding lemma we see that the subgroup $K_l\cap Ad_{\sigma(c,k)}K_k$ consists of matrices of the form
\begin{equation}
\left(\begin{array}{ccc}A_{11}&0&0\\
                          0&{A'}&0\\ 0&0&\mathbb J A_{11}\mathbb J
                               \end{array}\right)
\end{equation}
where
$$
A'=Ad_{\sigma(c,k-l)} \left(\begin{array}{ccc}A_{22}&0&0\\
                          0&B_{11}&B_{12}\\
                          0&B_{21}&B_{22}\end{array}\right)
$$
such that the whole determinant is $1$, with blocks $A_{11}\in
U(l)$, $A_{22}\in U(k-l)$, $B_{11}\in M_{n-2k}(\mathbb C)$, $B_{22}\in U(k-l)$ . When
$l=k$ we get matrices
\begin{equation}\label{Stief}
\left(\begin{array}{ccc}A_{11}&0&0\\
                          0&B_{11}&0\\
                       0&0&\mathbb JA_{11}\mathbb J
                               \end{array}\right)
\end{equation}
with $A_{11}\in U(k)$, $B_{11}\in U(n-2k)$ and $\det A_{11}^2=\det
B^{-1}$. In this case, applying lemma \ref{embed} exactly as in the
first part of the proof of Proposition \ref{nonstandardproj} we see that
$$
X_{l,k}=\frac{SU(n-l)}{\left\{\left(\begin{array}{cc}Ad_{\sigma(c,k-l) }A&0\\0&1_l\end{array}\right)\right\}}
$$
where $A\in S(U(k-l)\times U(n-k-l))$. Fix an auxiliary subgroup
$$
H_0=\left\{\left(\begin{array}{cc}\mathbf 1_l&0\\ 0&B\end{array}\right): B\in
 SU(n-l)\right\}
 $$
 and notice that $H_0/(H_0\cap Ad_{\sigma(c,k)} K_k)$ is Poisson diffeomorphic
to $K_l/(K_l\cap Ad_{\sigma(c,k)}K_k)$. Next $H_0/(H_0\cap Ad_{\sigma(c,k)} K_k)$ is easily seen to be Poisson
diffeomorphic to the standard Poisson quotients listed in the statement.
In the special case $l=k$ this yields the special case of Stiefel manifolds.
The symmetry provided by lemma \ref{symmetry} implies that the above results hold for $l\ge n-k$.

Now we consider the remaining case. Take $k<l<n-k$ (and hence $k\ne n/2$).
Then the intersection  $K_l\cap Ad_{\sigma(c,k)} K_k $ is given by:
\begin{equation}
\left\{\left(\begin{array}{cccc}A_{11}&0&0&0\\
                          0&B_{11}&0&0\\
                          0&0&B_{22}&0\\
                          0&0&0&{\mathbb J}_kA_{11}{\mathbb J}_k
                               \end{array}\right)\right\}
\end{equation}
such that the whole determinant is one and $A_{11}\in U(k)$, $B_{11}\in U(l-k)$, $B_{22}\in U(n-l-k)$.
We have, then,  considering that $K_l\cup J'$ generates $U(l)\times U(n-l)$
and then applying lemma 8.
$$
X_{l,k}\simeq \frac{U(l)\times U(n-l)}{J'}
$$
where $U(l)\times U(n-l)$ has the product Poisson structure (of standard Poisson $U(i)$'s) and
$J'$ consists of matrices
$$
{\left(\begin{array}{cccc}A_{11}&0&0&0\\ 0&B_{11}&0&0\\
0&0&B_{22}&0\\0&0&0&{\mathbb J}A_{11}{\mathbb J}\end{array}\right)}\,
$$
with no restrictions on determinants
(hence $J'$ is a Poisson--Lie subgroup of $U(l)\times U(n-l)$).
We remark that
$$
\frac{U(l)\times U(n-l)}{1\times U(l-k)\times U(n-k-l)\times 1}\simeq V_k^l\mathbb C\times V_{k}^{n-l}\mathbb C
$$
with the product of standard Poisson structures on Stiefel manifolds on the right.
It is just a quotient by a Poisson--Lie subgroup of $U(l)\times U(n-l)$.

The canonical embedding
$$
1\times U(l-k)\times U(n-k-l)\times 1\subseteq J
$$
of Poisson--Lie groups induces a   $U(l)\times U(n-l)$--equivariant, surjective Poisson map
$V_k^l\mathbb C\times V_{k}^{n-l}\mathbb C\twoheadrightarrow X_{l,k}$, with fibre $U(k)$.

Since the actions of the subgroups $1\times U(l-k)\times U(n-k-l)\times 1$ and $U(k)=\{A\oplus I_{k-1}\oplus
I_{n-k-1}\oplus A\, : A\in U(k)\}$ commute, $U(k)$ gives a well defined diagonal action on
$V_k^l\mathbb C\times V_{k}^{n-l}\mathbb C$
such that the quotient map onto its orbit space coincides with the above quotient map.

The codimension statement is can be easily verified by computation.
\fidi

\vspace{0.5cm}

\emph{\underline{Remarks}}
\begin{enumerate}
\item Explicit equations in Pl\"ucker coordinates for the embedded Poisson submanifolds $p_{\sigma_c}(K_l)$
can be obtained as in the $k=1$ case and are, at this point, matter of direct computations.
\item When $l<k$ or  $l>n-k$ the Poisson homogeneous spaces $X_{k,l}$ are the non standard version of
the Poisson homogeneous spaces denoted by $U/K^0_S$ in \cite{Stok}, where $U=SU(n)$ and $S$, subset of the
set of simple roots $\{\alpha_1,\ldots,\alpha_{n-l-1}\}$, in $\mathfrak{su}(n)$, is given
by deleting $\alpha_{k-l}$.
Such Poisson manifold should be compared with the standard Poisson quotient
$SU(n-l)/S(U(k-l)\times U(n-k-l))\times 1_l$ in the sense of understanding whether the two belong
to a Poisson pencil, as it is the case for projective and Grassmann manifolds.
\item As a last remark let us consider the maximal torus $T$ in
$SU(n)$ then $\pi\big|_T=0$. This implies that $T/(T\cap Ad_{\sigma(c,k)}
K_k)$ is a family of $0$--dimensional symplectic leaves in
Grassmannians which can be explicitly described:
$$
T\cap Ad_{\sigma(c,k)} K_k\simeq {\mathbb T}^{n-k}\Rightarrow
\frac{T}{T\cap Ad_{\sigma(c,k)} K_k}\simeq{\mathbb T}^k
$$
where the image of such points in the Grassmannian is given by
$$
(t_{1},\ldots,t_{k})\mapsto \langle t_{j}(\sqrt c e_{n-j+1,n-j+1}-
\sqrt{1-c}e_{j,n-j+1})\big| j=1,\ldots, k\rangle\, .
$$
\end{enumerate}

\end{document}